\documentclass[10pt]{article}          
\usepackage{graphicx}
\usepackage{url}
\usepackage{amssymb,amsthm}

\usepackage[caption=false]{subfig}
\usepackage{tikz,pgfplots}
\usepackage{stmaryrd}

\usepackage{bbm}
\usepackage{amsmath,amsfonts,bm} 
\usepackage{algorithm}
\usepackage{algpseudocode}

\usepackage{booktabs}

\usepackage{geometry}
\newcommand\lag{\mathcal{L}}
\newcommand\Hess{\mathcal{H}}
\newcommand\x{\mathbf{x}}

\DeclareMathOperator*{\argmin}{argmin}
\algdef{SE}[DOWHILE]{Do}{doWhile}{\algorithmicdo}[1]{\algorithmicwhile\ #1}%

\def\presuper#1#2%
  {\mathop{}%
   \mathopen{\vphantom{#2}}^{#1}%
   \kern-\scriptspace%
   #2}

\begin{document}

\title{$\mathbf{\frac 3 4}$-discrete optimal transport}
\author{ Fr\'ed\'eric de Gournay  \and Jonas Kahn \and L\'eo Lebrat}
\date{\today}
\maketitle
\begin{abstract}
This paper deals with the $\frac34$-discrete $2$-Wasserstein optimal transport between two measures, where one is supported by a set of segment and the other one is supported by a set of Dirac masses. We select the most suitable optimization procedure that computes the optimal transport and provide numerical examples of approximation of cloud data by segments.
\end{abstract}
%
%
%


\section*{Introduction}


The numerical computation of optimal transport in the sense of the 2-Wasserstein distance has known several breakthroughs in the last decade. One can distinguish three kinds of methods : The first method is based on underlying PDEs~\cite{benamou2000computational} and is only available when the measures are absolutely continuous with respect to the Lebesgue measure. The second method deals with discrete measures and is known as the Sinkhorn algorithm~\cite{cuturi2013sinkhorn,benamou2015iterative}. The main idea of this method is to add an entropic regularization in the Kantorovitch formulation. 
The third method is known as semi-discrete \cite{merigot2011multiscale,levy2015numerical,de2012blue} optimal transport and is limited to the case where one measure is atomic and the other is absolutely continuous. This method uses tools of computational geometry, the Laguerre tessellation which is an extension of the Voronoi diagram. The aim of this paper is to develop an efficient method to approximate a discrete measure (a point cloud) by a measure carried by a curve. This requires a novel way to compute an optimal transportation distance between two such measures, none of the methods quoted above  comply fully to this framework.

\vspace{2ex}

{\noindent \bf Pre-existing methods.}
The first method, the PDEs formulation, requires the measures to be absolutely continuous, which is not the case.

The Sinkhorn algorithm can be applied to this problem if the curve is sampled by points. But to the best of our knowledge \cite{Fred2017Gretsi}, taking Dirac masses along anisotropic ''objects'', here, a curve may dwindle the efficiency of this algorithm. Indeed, in this particular case the parameter of regularization of the Sinkhorn algorithm has to be chosen smaller than the curve-sampling precision, which causes numerical issues.

The semi-discrete optimal transport is the more favorable setting for such problem. In a previous paper \cite{Fred2018Projections} we described method to approximate an image by a curve. The method developed in \cite{Fred2018Projections} relies on a sampling of the curve by points and computes the Wasserstein distance between a density and a point cloud. This approach suffers from several drawbacks, the first one is that solving the optimal transport problem is more and more difficult as discretization step decreases. Indeed the closest the Dirac masses are taken along the curve the more stretched the Laguerre tessellation is. The second flaw is the poor quality of approximation of a curve by a sum of Dirac masses. The third problem is the need to transform the discrete measure into a measure with regular density.

Using ideas from the semi-discrete optimal transport, a solution to the above objections is to compute the optimal transport between an absolutely continuous measure and a measure carried by segments. The main tool of such computations is an extension of the Laguerre tessellation with conic boundaries. The difficulty arises when it comes to integrates the continuous density over theses Laguerre cells. In fact, a robust algorithm with exact integration is complex to develop.  
This is the reason why computer graphic community \cite{hoff1999fast,hiller2003beyond} implements numerical approximations for those Voronoi cells (Laguerre cells with equals weights) using shape primitives and rasterisation with graphic hardware. However numeric precision is intrinsically tied to grid (pixel). Hence, refinement for this method scales badly when the dimension increases.

Based on these observations we prefer to approximate the measure supported by curve by a measure supported by segments and to compute the optimal transport between an atomic measure and a measure supported on a set of segments. We coined this problem the $\frac 34 $-discrete optimal transport.

\vspace{2ex}

{\noindent \bf Contributions.}
This paper provides and studies an efficient algorithm to compute the 2-Wasserstein distance between a discrete measure and a measure supported by a set of segments. This algorithm is scalable in 2D and 3D and is parallelized. We also study the problem of optimizing the parameter of the measure carried by the set of segments and provide the formula of the gradient with respect to the parameter of the measure. 
The $\frac 3 4 $-discrete optimal transport benefits from the strength of Laguerre tessellation. Moreover it avoids the integration of the density over convex polygons (2D) or polyhedrons (3D), a common difficulty of semi-discrete optimal transport. Indeed the integrations needed for the computation of the cost function boils down to computing intersections between polyhedrons and to integrate the moments of $\nu$ over segments.

The ideas of this paper owes a much to  \cite{de2012blue,merigot2011multiscale,levy2015numerical} and the semi-discrete approach but it considers a measure supported by a set of segments instead of a measure with regular density. This modification annihilates the convergence theory developed in \cite{kitagawa2016convergence}. The optimal transport plan is no longer unique and Kantorovitch functional is no longer $\mathcal{C}^2$. But under slight geometric condition (see \eqref{eq::generic}) one can show that the dual functional is $\mathcal{C}^1$ with respect to the dual variables. 
\vspace{2ex}

{\noindent \bf Outline of the paper.}
Section \ref{sec::one} is devoted to set up the notations and the known results used in this paper. In Section \ref{sec::two} we compute the different derivatives of the cost function, in Section \ref{sec::three} we discuss the effective numerical implementation of the algorithm. In Section \ref{sec::four} we settle for the optimization procedure that solves the optimal transport, and in Section \ref{sec::five} we showcase numerical approximations of cloud data by segments. 

\section{Setting}
\label{sec::one}
\subsection{Reminders on 2-Wasserstein distance}
The 2-Wasserstein distance is a special instance of the optimal transport distance between two probabilities measures. It is defined as follow:

\paragraph{Given $\Omega \subset \mathbb{R}^d$, $\mu \in \mathcal{P}(\Omega)$ and $\nu \in \mathcal{P}(\Omega)$, the 2-Wasserstein distance between $\mu$ and $\nu$, $W_2(\nu,\mu)$ is given by:}
\begin{equation}\label{eq::MongeKanto} \tag{MK}
 W^2_2(\nu,\mu) = \inf_{\gamma \in \Pi(\nu,\mu)} \int_{\Omega^2} \|x-y\|^2_2 \quad \text{d} \gamma(x,y),
\end{equation}
where $\Pi(\nu,\mu)$ is the set of coupling between $\mu$ and $\nu$, that is, the set of measures whose marginals are $\mu$ and $\nu$ :
\begin{equation}\label{eq::coupling}
 \gamma \in \Pi(\nu,\mu) \Leftrightarrow \left\{ \begin{array}{l}
                                                  \displaystyle{\int_\Omega \psi(x) \text{d} \nu(x) = \int_{\Omega^2} \psi(x) \text{d} \gamma(x,y) \quad \forall \psi \in L^1(\nu)}\\
                                                  \text{and}\\
                                                  \displaystyle{\int_\Omega \phi(y) \text{d} \mu(y) = \int_{\Omega^2} \phi(y) \text{d} \gamma(x,y) \quad \forall \phi \in L^1(\mu)}
                                                 \end{array}\right.
\end{equation}

If both $(\Omega,\mu)$, $(\Omega,\nu)$ are Polish spaces, an elegant way \cite{villani2008optimal} to solve \eqref{eq::MongeKanto} is via its dual
\begin{align}\label{eq::Dual}
 \sup_{\psi \in L^1(\nu) ,\phi \in L^1(\mu)} & \int_\Omega \psi \text{d} \nu + \int_\Omega \phi \text{d} \mu \\
 \text{s.t} \quad \forall (x,y) \in \Omega^2& :  \quad \psi(x) + \phi(y) \leq \|x -y \|^2_2, \nonumber
\end{align}
where $\phi$ and $\psi$ are the Lagrange multipliers for \eqref{eq::coupling}, the marginals constraints of \eqref{eq::MongeKanto}.
Introducing the $c$-transform of $\phi$ as :
\[
\phi^c(x) = \inf_{y \in \Omega} \| x - y \|^2_2 - \phi(y).
\] 
The problem~\eqref{eq::Dual} can be rewritten as :
\begin{equation}\label{eq:tempoLaguerre}
 \sup_{\phi \in \mathbf{\Phi}^c(\Omega)} \int_\Omega \phi^c  \text{d} \nu + \int_\Omega \phi  \text{d} \mu ,
\end{equation}
where $\mathbf{\Phi}^c(\Omega)$ is the set of c-concave function on $\Omega$, see \cite{villani2008optimal}. Consider now the case where the measure $\mu$ is atomic~: $\mu(\x) = \sum_{i = 1}^n m_i \delta_{\x_i}, \x_i \in \mathbb{R}^d, m_i \in \mathbb{R}$. In this case $\phi$ belongs to $\mathbb{R}^ n$, and $\phi^c(x) = \min_{i \in \llbracket 1 , n\rrbracket} \|x - \x_i \|^2_2 - \phi_i $.
This naturally leads to the definition of the $i$-th Laguerre cell~\cite{aurenhammer1987power} :
\begin{equation}
\label{eq::Laguerre}
\lag_i(\x,\phi)=\{x\in \Omega \text { such that } \| x - \x_i\|^2_2- \phi_i \le \|x- \x_j\|^2_2 - \phi_j \quad \forall j\in \llbracket 1,n \rrbracket \}.
\end{equation}
%
Provided that $\nu(\lag_i \cap \lag_j) = 0$ for every $i\ne j$, the final problem states as :

\begin{equation} \label{eq::costfunc}
 W^2_2(\mu,\nu) = \sup_{\phi \in \mathbb{R}^n} g(\phi,\x), \text{ with } g(\phi,\x) = \sum_i \int_{\lag_i (\x,\phi)} \left( \|x - \x_i\|^2_2 - \phi_i \right) d\nu + \sum_i m_i \phi_i.
\end{equation}
Denoting $\phi^\star$ as the solution of maximization problem \eqref{eq::costfunc} the physical interpretation of the Laguerre cell $\lag_i(\x,\phi^\star)$  is that the Dirac mass located at $\x_i$ is transported to $\text{supp}(\nu)\cap \lag_i$.
There is equality between problem \eqref{eq::MongeKanto} and problem \eqref{eq::costfunc} since the quadratic cost is continuous and $(\Omega,\nu)$,$(\Omega,\mu)$ are two Polish spaces \cite{villani2008optimal}.

\subsection{Setting}
In this paper $\mu$ will denote a $n$-atomic probability measure on $\mathbb{R}^d$ :
\[
\mu(\x) = \sum_{i = 1}^n m_i \delta_{\x_i}, \text{ s.t} \quad \sum_i m_i = 1 \text{ with} \quad \x_i \in \mathbb{R}^d, m_i \in \mathbb{R}^+, \text{ for all } 1 \leq i \leq n.
\]
By contrast, $\nu$ will be a probability measure supported by a polyline,  
\begin{equation} \label{eq::polyLine}
 \nu = \sum_{\alpha=1}^{p}  \rho_\alpha l^\alpha_{\#}  \lambda_{[0,1]}  \quad \text{s.t.} \quad
 \begin{array}{rcll}
 &&&\\
 l^\alpha:& [0,1] &\rightarrow& \mathbb{R}^d\\
 &t  &\mapsto& (1-t) P_\alpha + t P_{\alpha+1}                                                                                                                                                                      
 \end{array}\quad , \rho_\alpha \in \mathbb{R}^{+}, P_\alpha \in \mathbb{R}^d,
\end{equation}
where the measure $ l^\alpha_{\#}  \lambda_{[0,1]} $ is the push-forward through $l^\alpha$ of $\lambda_{[0,1]}$, the $[0,1]-$Lebesgue measure. 
The measure $l^\alpha_{\#}  \lambda_{[0,1]}$ is defined as :
\[
l^\alpha_{\#}  \lambda_{[0,1]} (A) := \lambda_{[0,1]} \left( (l^\alpha) ^{-1} (A) \right) \text{ for every Borelian set } A. 
\]

The fact that $\nu$ is a probability measure translates into $\sum_\alpha \rho_\alpha = 1$.
Notice that $\rho_\alpha=0$ is admissible and in this case the support of the resulting probability $\nu$ will be composed of disjoint polylines.

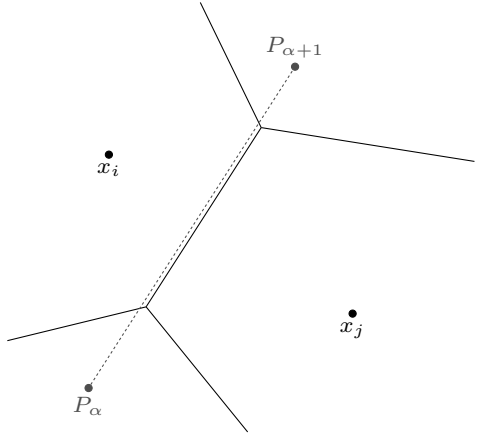
\begin{figure}[!ht]
 \centering
 \begin{tikzpicture}[scale=4.5]
\definecolor{xdxdff}{rgb}{0.49,0.49,1}
\definecolor{qqqqff}{rgb}{0,0,1}
\clip(0,0.11) rectangle (1.45,1.43);
\draw (0.78,1.03)-- (0.44,0.5);
\draw (0.6,1.4)-- (0.78,1.03);
\draw (0.78,1.03)-- (1.41,0.93);
\draw (0.44,0.5)-- (0.03,0.4);
\draw (0.44,0.5)-- (0.74,0.13);
\draw[color=black!70,dash pattern=on 1pt off 1pt] (0.88,1.21)-- (0.27,0.26);
\begin{scriptsize}
\draw[fill] (0.33,0.95) circle   (.3pt) node[below,scale=1.2] {$x_i$};
\draw[fill] (1.05,0.48) circle   (.3pt) node[below,scale=1.2] {$x_j$};
\draw[fill,color=black!70] (0.27,0.26) circle   (.3pt) node[below,scale=1.2] {$P_\alpha$};
\draw[fill,color=black!70] (0.88,1.21) circle   (.3pt) node[above,scale=1.2] {$P_{\alpha+1}$};
\end{scriptsize}

\end{tikzpicture}
 \caption{When the hypothesis \eqref{eq::generic} is violated one can cook up a $\phi$ such that the functional $g$ in \eqref{eq::costfunc} is no longer $\mathcal{C}^1$ with respect to $\phi$. Here a small perturbation on $\phi_i$ or $\phi_j$ induces a discontinuity on $\nabla_\phi g$.}
 \label{fig::nongeneric}
\end{figure}

In order to ensure the regularity of function $g$ defined in \eqref{eq::costfunc} we enforce hypothesis \eqref{eq::generic} throughout the paper.


\begin{equation} \label{eq::generic}
 \forall \alpha \in \llbracket 1, p \rrbracket , \quad \forall (i,j) \in \llbracket 1, n \rrbracket^2, \quad  i \neq j : \quad \langle P_{\alpha+1} - P_\alpha , \x_i - \x_j \rangle \neq 0. \tag{H}
\end{equation}
Under \eqref{eq::generic}, for all $\phi$ we have $\nu(\lag_i(\x,\phi) \cap \lag_j(\x,\phi)) = 0 $. Indeed if $\lag_i$ and $\lag_j$ are two Laguerre cells with a common boundary $\lag_i \cap \lag_j$, the boundary has to be orthogonal to $\x_i - \x_j$. Hypothesis \eqref{eq::generic} prevents tiny perturbation of $\phi$ to harshly shift the affectation of the segment's mass from one Dirac mass to another, see Figure~\ref{fig::nongeneric}. Under this assumption $g$ is a $\mathcal{C}^1$ function of $\phi$, see \cite{Fred2018differentiation}.


Let $\mathcal I$ be the set of indices $(i,\alpha)$ such that the $i^\text{th}$ Laguerre cell meets the $\alpha^\text{th}$ segment, that is: 
\[\mathcal{I} = \left\lbrace (i,\alpha) \text{ s.t. } \lag_i(\x,\phi) \cap [P_\alpha,P_{\alpha+1} ] \neq \emptyset \right\rbrace.\] If the space $\Omega$ is convex then the Laguerre cells defined in \eqref{eq::Laguerre} are convex and $\lag_i \cap [P_\alpha,P_{\alpha+1}] $ is a segment. For every $(i,\alpha) \in \mathcal{I} $ denote $t_s^{i\alpha},t_e^{i\alpha}$ the starting (respectively ending) time of the $\alpha^{\text{th}}$ segment $[P_\alpha,P_{\alpha+1}]$ in the $i^\text{th}$ Laguerre cell :
\begin{equation}\label{eq::intersections}
 0 \leq t^{i\alpha}_s \leq t^{i\alpha}_e \leq 1 \text{ and }\lag_i \cap [P_\alpha,P_{\alpha+1}] = [l^\alpha(t^{i\alpha}_s), l^\alpha(t^{i\alpha}_e)].
\end{equation}
Then the cost function $g(\phi,\x)$ defined in \eqref{eq::costfunc} can be re-written as :
%
\begin{align}
 g(\phi,\x) &= \sum_{(i,\alpha) \in \mathcal{I}}   \int_{t^{i\alpha}_s}^{t^{i\alpha}_e} r^{i\alpha}(t) dt  + \sum_{i = 1}^n \phi_i m_i  \label{eq::costFunctionComputer} \\
 \text{ with } r^{i\alpha}(t) &= (\| \l^\alpha(t) - \x_i \|^2_2 - \phi_i ) \rho_\alpha \nonumber 
\end{align}

\section{Derivatives of the cost function}
\label{sec::two}
\subsection{First order derivative with respect to $\phi$}
Hypothesis \eqref{eq::generic} ensures that $\nu(\lag_i(\x,\phi)\cap \lag_j(\x,\phi)) = 0, \forall i\neq j.$ A direct application of \cite{Fred2018differentiation} shows that $g$ is differentiable with respect to $\phi$ and : 
\begin{equation} \label{eq::FirstOrderDerivative}
\frac{\partial g }{\partial \phi_i} (\phi,\x) = m_i - \int_{\lag_i(\x,\phi)} d\nu.
\end{equation}

\subsection{Computation of the second order derivative with respect to $\phi$}
\label{sec::Hess}
If there exists a segment $[P_\alpha,P_{\alpha+1}]$ that passes through the intersection of at least three Laguerre cells then the functional $g$ fails to be twice differentiable. On the one hand it is impossible to design an hypothesis in the spirit of hypothesis \eqref{eq::generic} that can prevent such a pathological case to happen during the optimization in $\phi$. On the other hand such a baneful case almost surely never happens. Hence in this section the Hessian is computed without proof of existence by following a cumbersome calculus.

\begin{figure}[ht]
\begin{minipage}[b]{0.35\linewidth}
\begin{tikzpicture}[xscale=2.5,yscale=2.9]
\coordinate (I) at (0,0);
\coordinate (G) at ( .74666667,.84);
\coordinate (H) at (1.25333334,.84);
\coordinate (J) at (2.,0);
\coordinate (O) at (1.,.8);
\coordinate (P) at (1.,.84);
\coordinate (Q) at (0.623448275862069,0.601379310344828);
\coordinate (R) at (0.573793103448276,0.64551724137931);
\coordinate (S) at (1.376692057130982,0.601362793517391);
\coordinate (T) at (1.426330170161779,0.645516960094837);

\draw[fill] (1,.4) circle   (.3pt) node[below,scale=.8] {$x_j$};
\draw[fill] (.4,.8) circle  (.3pt) node[left,scale=.8] {$x_l$};
\draw[fill] (1,1.5) circle  (.3pt) node[below,scale=.8] {$x_k$};
\draw[fill] (1.6,.8) circle (.3pt) node[right,scale=.8] {$x_m$};
\draw (0,-.1) -- (.8,.8) -- (0,1.4);
\draw (.8,.8) -- (1.2,.8);
\draw (2,-.1) -- (1.2,.8) -- (2,1.4);

\draw[color=gray!70] (I) -- (G) -- (H) -- (J);
\draw[->,thick,color=gray!90] (O) -- (P) node[yshift=-0.4cm,scale=.9] {$\delta_{jk}$};
\draw[->,thick,color=gray!90] (Q) -- (R) node[xshift=.3cm,yshift=-0.4cm,scale=.9] {$\delta_{jl}$};
\draw[->,thick,color=gray!90] (S) -- (T) node[xshift=-.27cm,yshift=-0.4cm,scale=.9] {$\delta_{jm}$};

\draw [color=black!80](1.8,.6) node[scale=.8] {$\mathcal{L}_m$};
\draw [color=black!80](1.020189314527146,0.135114508213524) node[scale=.8] {$\mathcal{L}_j$};
\draw [color=black!80](1.020189314527146,1.185114508213524) node[scale=.8] {$\mathcal{L}_k$};
\draw [color=black!80](0.25,.6) node[scale=.8] {$\mathcal{L}_l$};

\end{tikzpicture}

\end{minipage}
\hspace{1cm}
\begin{minipage}[b]{0.55\linewidth}
\centering
 \begin{tikzpicture}[xscale=10.,yscale=20]
\definecolor{polyCol}{rgb}{.4,0.4,0.4}
\coordinate (E) at (0.8,0.8);
\coordinate (F) at (1.2,0.8);
\coordinate (G) at (0.746666666666667,0.84);
\coordinate (H) at (1.253333333333333,0.84);
\coordinate (V) at (0.666666666666667,0.9);
\coordinate (W) at (1.333333333333333,0.9);
\coordinate (A1) at (0.622222222222222,0.7);
\coordinate (B1) at (0.711111111111111,0.7);
\coordinate (C1) at (1.288952145602416,0.7);
\coordinate (D1) at (1.377866337176715,0.7);

\coordinate (E1) at (0.853442657813992,0.875292532674489);

\coordinate (F1) at (1.105440770132722,0.731114508213524);
\coordinate (G1) at (1.434213677660114,0.731114508213524);

\coordinate (H1) at (0.98281899051799,0.8);
\coordinate (I1) at (0.98281899051799,0.84);

\coordinate (K1) at (0.928649426082078,0.84);

\coordinate (P1) at (1.015681899051799,0.83);

\coordinate (J1) at (1.261275122952846,0.731114508213524);
\coordinate (L1) at (1.335913235983642,0.74826867479097);
\coordinate (M1) at (1.350189314527146,0.731114508213524);
\coordinate (P2) at (1.310189314527146,0.745114508213524);

\coordinate (J2) at (1.283,0.706114508213524);
\coordinate (M2) at (1.373,0.706114508213524);

\coordinate (HH) at (1.13581899051799,0.8);
\coordinate (II) at (1.13581899051799,0.84);


\draw (E)--(F);
\draw (V)--(E)--(B1);
\draw (W)--(F)--(C1);
\draw[color=gray!70] (A1) -- (G) -- (H) -- (D1);
\draw[dashed,color=polyCol] (E1) -- (F1) -- (G1);
\draw[fill,color=black!70] (E1) ellipse (.002cm and .001cm) node[above,scale=0.8] {$P_{\alpha+1}$};
\draw[fill,color=black!70] (F1) ellipse (.002cm and .001cm) node[below,scale=0.8] {$P_{\alpha}$};
\draw[fill,color=black!70] (G1) ellipse (.002cm and .001cm) node[below,scale=0.8] {$P_{\alpha-1}$};

\draw[->] (H1)--(I1) node[xshift=.25cm,yshift=-0.35cm,scale=.9] {$n_{jk}$};

\draw[<->] (HH)--(II)node[xshift=-.22cm,yshift=-0.35cm,scale=.9] {$\delta_{jk}$};

\draw[<->] (J2)--(M2)node[xshift=-.37cm,yshift=-0.25cm,scale=.8] {$\frac{\partial g}{\partial \phi_j \partial \phi_m}$};

\coordinate (AA) at (0.97281899051799,0.8);
\coordinate (BB) at (0.90281899051799,0.84);
\draw[<->] (AA)--(BB)node[xshift=-.22cm,yshift=-0.50cm,scale=.8] {$\frac{\partial g}{\partial \phi_j \partial \phi_k}$};

\begin{scope}
 \path[clip] (H1)--(I1)--(K1);
 \fill[gray,opacity=.6,draw=black] (H1) circle (0.5pt);
\end{scope}
\draw [black,opacity=.6](P1) node[scale=.9,xshift=-.65cm,yshift=-.00cm] {$\theta^{jk}_{\alpha}$};

\begin{scope}
 \path[clip] (J1)--(L1)--(M1);
 \fill[gray,opacity=.6,draw=black] (J1) circle (0.7pt);
\end{scope}

\coordinate (JJ) at (1.217275122952846,0.781114508213524);
\coordinate (LL) at (1.2900913235983642,0.798268675);
\draw[<->] (JJ)--(LL) node[xshift=-.4cm,yshift=-.00cm,scale=.9,rotate=20] {$\delta_{jm}$};

\draw[->] (J1)--(L1) node[xshift=-.43cm,yshift=-.050cm,scale=.9,rotate=20] {$n_{jm}$};
\draw [black,opacity=.6](P2) node[scale=.9,xshift=.10cm,yshift=-.57cm] {$\theta^{jm}_{k-1}$};

\draw [color=black!80](1.360189314527146,0.795114508213524) node[scale=1.] {$\mathcal{L}_m$};
\draw [color=black!80](0.920189314527146,0.735114508213524) node[scale=1.] {$\mathcal{L}_j$};
\draw [color=black!80](1.180189314527146,0.885114508213524) node[scale=1.] {$\mathcal{L}_k$};
\end{tikzpicture}
\end{minipage}

\caption{(Left) Normal displacements $(\delta_{il},\delta_{ik},\delta_{im})$ of the boundary of the Laguerre cell when $\phi_i$ increases and $(\phi_l,\phi_k,\phi_m)$ remain constant. (Right) geometrical interpretation of the Hessian}
\label{fig::displacementBoundary}
\end{figure}
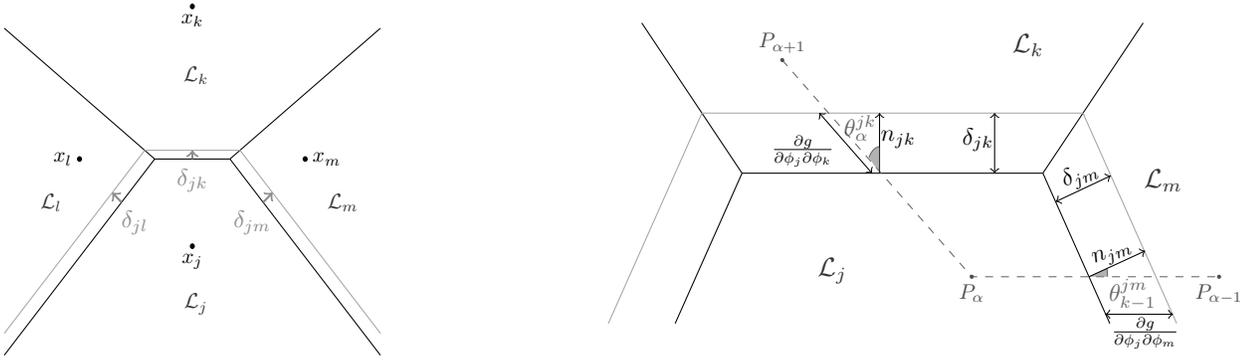

Denote $(e_i)_i$ the canonical basis. Following the computation of the first derivative given in \eqref{eq::FirstOrderDerivative}, the second order derivative is given by :
\begin{equation}
 \frac{\partial^2 g}{\partial \phi_i \phi_j} = -  \lim_{\varepsilon_i \rightarrow 0} \frac{\int_{\lag_i(\x,\phi + \varepsilon e_j)} d\nu - \int_{\lag_i(\x,\phi) } d\nu}{\varepsilon}.
\end{equation}
Denote $n_{il}$ the outer normal of $\lag_i$ on the facet $\lag_i\cap \lag_l$, that is $n_{il} = \frac{ x_l -x_i}{\| x_l -x_i \|}$. Denote by $\delta_{il}$ the first order approximation of the evolution of the facet $\lag_i \cap \lag_l$ in the normal direction $n_{il}$ when we change the Lagrange multiplier of the $j$-th cell,  that is, when we change $\phi$ into $\phi + \varepsilon e_j$, see Figure \ref{fig::displacementBoundary} (left). If both $i$ and $l$ are different from $j$ then $\delta_{il} = 0$ otherwise $\delta_{ij}$ is given by :
\begin{align}
 &y + \varepsilon \delta_{ij} n_{ij}  + o(\varepsilon) \in \left(\lag_i \cap \lag_j  \right) (\x,\phi + \varepsilon e_j) \nonumber \\ 
 \Leftrightarrow& \| y + \varepsilon \delta_{ij} n_{ij} -\x_j \|_2^2 - \phi_j - \varepsilon  + o(\varepsilon) =  \| y + \varepsilon \delta_{ij} n_{ij} -\x_i \|_2^2 -\phi_i\label{eq::OnEdge} \\
 \Rightarrow&\delta_{ij} = \frac{1}{2 \| \x_j -\x_i \|} . \nonumber
\end{align}

Denote $\Theta_{ij}$ the indices of the segments that intersect the facet $\lag_i \cap \lag_j$, it is given by 
\[
\Theta_{ij} = \left\lbrace \alpha \text{ such that, } (i,\alpha) \in \mathcal I \text{ and } (j,\alpha) \in \mathcal I \right\rbrace. 
\]
For every $\alpha$ in $\Theta_{ij}$ denote $\theta_\alpha^{ij}$ the angle between the segment $[P_\alpha,P_{\alpha+1}]$ and the outer normal $n_{ij}$. See Figure~\ref{fig::displacementBoundary} (right), simple geometrical consideration shows that
the Hessian formula boils down to :
\begin{equation}\label{eq::Hessian}
 \frac{\partial g}{\partial \phi_i \phi_j}(\phi,\x) = \left\{ \begin{array}{l|l}
		     \displaystyle{  - \sum_{\theta \in \Theta_{ij}}  \frac{\rho_\alpha}{2 \|\x_j - \x_i \| \cos(\theta_\alpha^{ij}) }} & \text{if} \ i \neq j ,\\               
                     \displaystyle{ - \sum_{k \neq i} \frac{\partial g}{\partial \phi_i \phi_k}  }& \text{otherwise.}
                                                 \end{array}\right.
\end{equation}

%

\subsection{Computation of the first order derivative with respect to $P$ and $\rho$}

In this section we compute the derivatives of the 2-Wasserstein distance with respect to the parameters of the measure $\nu$. We denote generically as $\partial_\nu$ the derivative with respect to one of the parameters that define $\nu$, that is either one of the positions $P_\alpha$ or one of the densities $\rho_\alpha$. Let $\phi^\star$ be the optimal Lagrange multiplier in \eqref{eq::costfunc}
\[
 G(\nu) := g(\phi^\star,\nu)=W_2^2(\mu,\nu)
\]

The derivative of $g(\phi^\star,\nu)$ with respect to the parameters of $\nu$ are given using the chain rule formula.
\begin{equation}\label{eq::DerivNu}
 \frac{\text{d}G}{\text{d}\nu} = \frac{\partial g}{\partial \phi} \frac{\partial \phi}{\partial \nu} + \frac{\partial g}{\partial \nu}
\end{equation}
Since $\phi^\star$ is a solution of \eqref{eq::costfunc}, the derivative of $g$ with respect to $\phi$  is zero at $\phi^\star$. Hence computing the differential of $G$ consists in differentiating the cost function $g$  while keeping $\phi$ fixed at $\phi^\star$.
Differentiating \eqref{eq::costFunctionComputer} with respect to $\nu$ we obtain :
\[
 \partial_\nu G = \sum_{(i,\alpha) \in \mathcal{I}} (\partial_\nu t^{i\alpha}_e) r^{i\alpha}(t^{i\alpha}_e) -  (\partial_\nu t^{i\alpha}_s) r^{i\alpha}(t^{i\alpha}_s) +    \int_{t^{i\alpha}_s}^{t^{i\alpha}_e} \partial_{\nu} r^{i\alpha}(t) dt
\]

When considering $\partial_\nu t^{i\alpha}$ only three cases can occur :
\begin{enumerate}
 \item If $P_\alpha$ belongs to the interior of the ith Laguerre cell then, $t_s^{i\alpha} = 0$ and $\partial_\nu t_s^{i\alpha} = 0$.
 \item If $P_{\alpha+1}$ belongs to the interior of the ith Laguerre cell then, $t_e^{i\alpha} = 1$ and $\partial_\nu t_e^{i\alpha} = 0$.
 \item For all $i$ there exists exactly one $j$ such that $t_s^{i\alpha} =t_e^{j\alpha}$. The segment $[P_{\alpha},P_{\alpha+1}]$ intersects $\lag_i \cap \lag_j$ at point $l^{\alpha}(t_s^{i\alpha})$. For all $x$ in $\lag_i \cap \lag_j$, by definition of the Laguerre cell \eqref{eq::Laguerre}, the following equality holds $\| x - \x_i\|^2_2- \phi_i = \|x- \x_j\|^2_2 - \phi_j$ then $r^{i\alpha}(t^{i\alpha}_e) = r^{j\alpha}(t^{j\alpha}_s)$.
\end{enumerate}

It follows that 
\begin{equation}\label{eq::derivNu}
 \partial_\nu G =  \sum_{(i,\alpha) \in \mathcal{I}} \int_{t^{i\alpha}_s}^{t^{i\alpha}_e} \partial_{\nu} r^{i\alpha}(t) dt
\end{equation}


\subsubsection{Derivative with respect to $P$}
\label{sec::derivP}
We are first interested in differentiating the Wasserstein distance with respect to the position of the endpoints $(P_\alpha)_{\alpha \in \llbracket 1, p \rrbracket}$ of the polyline. A direct application of \eqref{eq::derivNu} yields :

\[
\frac{\partial G}{\partial P_\alpha} = \int_{t^{i(\alpha-1)}_s}^{t^{i(\alpha-1)}_e} \partial_{P_\alpha} r^{i(\alpha-1)}(t) dt + \int_{t^{i\alpha}_s}^{t^{i\alpha}_e} \partial_{P_\alpha} r^{i\alpha}(t) dt \\
\]

The differential of $r^{i\alpha}(t)$ with respect to $P_\alpha$ amounts to differentiate $\|l^\alpha(t) -x_i \|^2$, we recall that :
\[
\left\{ \begin{array}{ll}
      \displaystyle{\frac{\partial \| l^{\alpha}(t) -x_i\|^2}{\partial P_\alpha}} &= 2(1-t) \left( l^{\alpha}(t) -x_i \right)\\
      \text{}& \\
      \displaystyle{\frac{\partial \| l^{(\alpha-1)}(t) -x_i\|^2}{\partial P_\alpha}} &= 2t \left( l^{(\alpha-1)}(t) -x_i \right).
	\end{array}
\right.
\]
Let $\mathcal{K}(\alpha)$ the set of indices $i$ such that the $i$-th Laguerre cell has a non-empty intersection with the $\alpha$-th segment:
\[
 \mathcal{K}(\alpha) = \{ i \text{ s.t. } \lag_i(x,\phi) \cap [ P_\alpha,P_{\alpha+1} ] \neq \emptyset  \} = \{ i \text{ such that } (i,\alpha) \in \mathcal{I} \}.
\]
A straightforward computation yields

%
%

\begin{align}
\frac{\partial G}{\partial P_\alpha} =& \left( P_\alpha + \frac{P_{\alpha - 1} - P_\alpha}{3} - \sum_{i 
\in \mathcal{K}(\alpha-1)} \int^{t_e^{i(\alpha-1)}}_{t_s^{i(\alpha-1)}}  t \x_i dt\right) \rho_{\alpha-1} \nonumber\\
+& \left(P_\alpha + \frac{P_{\alpha + 1} - P_\alpha}{3} - \sum_{i \in \mathcal{K}(\alpha)} \int^{\presuper{\dag}{t_e^{i\alpha}}}_{\presuper{\dag}{t_s^{i\alpha}} } t\x_i dt\right)\rho_{\alpha} \label{eq::torque}
\end{align}

where $\presuper{\dag}{t} $ is the reverse parameterization of the segment $[P_{\alpha},P_{\alpha+1}]$, that is
\[
 \presuper{\dag}{t_s^{i\alpha}} = 1 - t_e^{i\alpha} \quad \text{and} \quad \presuper{\dag}{t_e^{i\alpha}} = 1 - t_s^{i\alpha}.
\]

We may interpret \eqref{eq::torque} as the sum of the torques of the surrounding segments $[P_{\alpha-1},P_{\alpha}]$ and $[P_{\alpha},P_{\alpha+1}]$ around the point $P_\alpha$.

\subsubsection{Derivative with respect to $\rho$}
\label{sec::derviRho}
The functional $g$ ~\eqref{eq::costFunctionComputer} is linear in $\rho_\alpha$, and its derivative is therefore trivial. Throughout this paper, we consider that the density is constant by segment. As a consequence 
\begin{equation}
\label{eq::ohputainyapasdelabel}
 \rho_\alpha= \frac{\| P_\alpha - P_{\alpha+1} \|}{\sum_{\alpha} \| P_\alpha - P_{\alpha+1} \|},
\end{equation}
hence the derivative of $\rho$ with respect to $P$ is easily computed.

\subsubsection{Optimization algorithm}
The aim of this section is to describe the algorithm which optimizes the Wasserstein distance $G$ with respect to the position of the endpoints $P_\alpha$ of the polyline.
The proposed algorithm is a gradient descent method with a metric $\Sigma$, we recall that the gradient of $G$ with respect to $P$ is given by :
\begin{equation} \label{eq::nablaG}
 \nabla_P G = \frac{\partial G}{\partial P} + \frac{\partial G}{\partial \rho} \frac{\partial \rho}{\partial P},
\end{equation}

The differential of $G$ with respect to $P$ and $\rho$ are discussed in Section \ref{sec::derivP} and Section \ref{sec::derviRho} respectively. The main goal of this section is to discuss the choice of the metric $\Sigma$.

Let us consider an isolated segment $[P_\alpha, P_{\alpha+1}]$, that is $\rho_{\alpha - 1} = \rho_{\alpha + 1} = 0$. The average of the gradient for the segment $[P_\alpha,P_{\alpha+1}]$ is given by :
\begin{equation} \label{eq::barycenterSimple}
 \frac 1 2 \left(  \frac{\partial G}{\partial P_\alpha}  + \frac{\partial G}{\partial P_{\alpha+1}}\right) = \left(\frac 1 2 P_\alpha + \frac 1 2 P_{\alpha +1 } - \sum_{i \in \mathcal{K}_\alpha} \int ^{t_e^{i\alpha}}_{t_s^{i\alpha}} \x_i  dt \right) \rho_\alpha = \rho_\alpha ( c_\alpha - \bar{\x}_\alpha)
\end{equation}

where $c_{\alpha}$ is the center of the $\alpha$-th segment and $\bar{\x}_\alpha$ is the weighted average of points $\x_i$ seen by the $\alpha$-th segment :

\[
  c_\alpha = \frac{P_{\alpha} + P_{\alpha + 1}}{2} \quad \text{and} \quad \bar{\x}_\alpha = \sum_{i \in \mathcal{K}_\alpha} \int_{t_s^{i\alpha}}^{t_e^{i\alpha}} \x_i \text{d}t.
\]
Equation \eqref{eq::barycenterSimple} is reminiscent of the formula of the derivative of $G$ in the semi-discrete setting, see \cite{levy2015numerical,de2012blue,merigot2011multiscale}. In this setting, $\mu$ is a measure with density and the approximating measure $\nu$ is a sum of Diracs~:
\[\nu=\sum_{\alpha = 1}^p m_\alpha \delta_{\bf{y}_\alpha}\]

In this case, the derivative of $G$ with respect to the position of the points is given by : 
\begin{equation}\label{eq::BlueNoise}
 \frac{\partial G}{\partial \bf{y}_\alpha} = m_\alpha(\bf{y}_\alpha - b_\alpha),
\end{equation}
where ${\bf b}_\alpha= \frac{1}{m_\alpha}\int_{\lag_\alpha} x d\mu$ is the barycenter of the $\alpha$th Laguerre cell.

In the semi-discrete setting, the most commonly used algorithm when minimizing $G$ with respect to $\bf{y}$ is to update the point position $\bf{y}_\alpha$ to their barycenters ${\bf b}_\alpha$. This procedure is known as Lloyd's algorithm \cite{de2012blue,merigot2011multiscale,levy2015numerical}. In view of the formula of the gradient~\eqref{eq::BlueNoise}, this method is a gradient descent with metric $\Sigma = \text{diag}(m_\alpha)$. In our case it is natural to consider a metric $\Sigma$ defined via $\rho_\alpha$ as an analogy to the semi-discrete setting. We define :
\begin{equation} \label{eq::metric}
\Sigma = \text{diag}(\frac{\rho_{\alpha-1} + \rho_\alpha}{2}). 
\end{equation}

The main algorithm that minimizes $G$ with respect to $P$ is given by algorithm \ref{alg::GradientDescent} below 
\begin{algorithm}
\caption{Optimization polyline position}\label{alg::GradientDescent}
\begin{algorithmic}[1]
\Procedure{Optimization in P}{}
\State $\nabla_P G \gets \mathbf{0}$
\Do
  \State $P \gets P -\Sigma^{-1} \cdot \nabla_P G $ \Comment{ with $\Sigma$ defined in \eqref{eq::metric}, $\nabla_P G$ defined in \eqref{eq::nablaG} }
  \State $\rho \gets \texttt{normalizedDensity}(P)$ \Comment{ as in \eqref{eq::ohputainyapasdelabel} }
  \State $\phi^*,\nabla_P G \gets \texttt{computeOptimalTransport}(\x,m,P,\rho)$ \Comment{ see Algorithm \ref{alg::GradientPhiDescent} }
\doWhile {$\| \nabla P \|_\infty \geq  10^{-3}$}
\EndProcedure
\end{algorithmic}
\end{algorithm}

\section{Numerical implementation}
\label{sec::three}
The Laguerre cells are computed using the computational geometry library CGAL ~\cite{cgal}. In $2$D the algorithm is fast and scalable, the average complexity for $n$ sites randomly drawn is linear in time and memory. For $3$D triangulation the worst case complexity is quadratic, but for random point configurations the complexity is observed to be  almost linear~\cite{dwyer1991higher}. In this section we discuss the computation of $g$ and its parallelization.

\subsection{Integration computation}
\label{sec::integrals}
The main issue when computing $g$ and its derivative is to evaluate the intersections~\eqref{eq::intersections} between the polylines and the Laguerre cells. The intersections are computed by following each segment $l^\alpha(t)$ with $t$ increasing. We first focus our attention on computing the exit time of the Laguerre cell $j$ knowing the starting time $t_s^{i\alpha}$. To clarify things suppose that at the known time $t_e^{i\alpha} = t_s^{j\alpha}$ the segment exits the Laguerre cell $\lag_i$ and enters the Laguerre cell $\lag_j$. The objective is to compute the time $t_e^{j\alpha} = t_s^{k\alpha}$ and the index $k$ such that the segment exits the Laguerre cell $\lag_j$ at time $t_e^{j\alpha}$ and enters the Laguerre cell $\lag_k$. Such a time is computed by solving the following minimization problem

\begin{equation}\label{eq::nextTime}
k = \argmin_{m \text{ s.t. }m\ne i, \  t_s^{j\alpha}<t_m \le 1} t_m, \quad \text{with} \quad t_m = \frac{2 \langle P_\alpha , \mathbf{x}_j - \mathbf{x}_m \rangle + \| \mathbf{x}_m\|^2 - \phi_m + \phi_j - \|\mathbf{x}_j\|^2_2}{2\langle P_{\alpha} - P_{\alpha+1} , \mathbf{x}_j - \mathbf{x}_m \rangle}. 
\end{equation}

The exit time $t_e^{j\alpha}$ is now equal to $t_k$. Note that the choice of $m$ can be restricted to the indexes of adjacent Laguerre cells, which is a small set in practice (a dozen of indexes).
If the set $\{ m\ne i, \  t_s^{j\alpha}<t_m \le 1 \}$ is empty, the segment ends in the Laguerre cell $\lag_j$ and we set $t_e^{j\alpha} = 1$ and we stop the procedure for the segment.

In the case when $t_s^{j\alpha}= 0$, that is the segment starts in the Laguerre cell $\lag_j$, we consider the same minimization problem as \eqref{eq::nextTime} where the constraint set is replaced with $\{m, \ 0 \leq t_m \leq 1\}$ .

We still have to compute the index of the Laguerre cell where the segment $[P_\alpha,P_{\alpha +1}]$ begins. If $\alpha > 0$, it is obviously the index of the ending Laguerre cell of $[P_{\alpha-1},P_{\alpha}]$. For the case $\alpha=0$, we add a dummy segment $[\x_l,P_\alpha]$, where $\x_l$ is the position of the $l$-th Dirac mass corresponding to the largest multiplier $\phi_l = \max_k \phi_k$. By definition of the Laguerre diagram ~\eqref{eq::Laguerre} the point $\x_l$ belongs to $\lag_l$.

Given the starting and ending times, the computation of $g(\phi,\x)$ and its derivative abridge to integrate polynomials within $[t_s^{j\alpha},t_e^{j\alpha}]$ using Gaussian quadrature.

\medskip

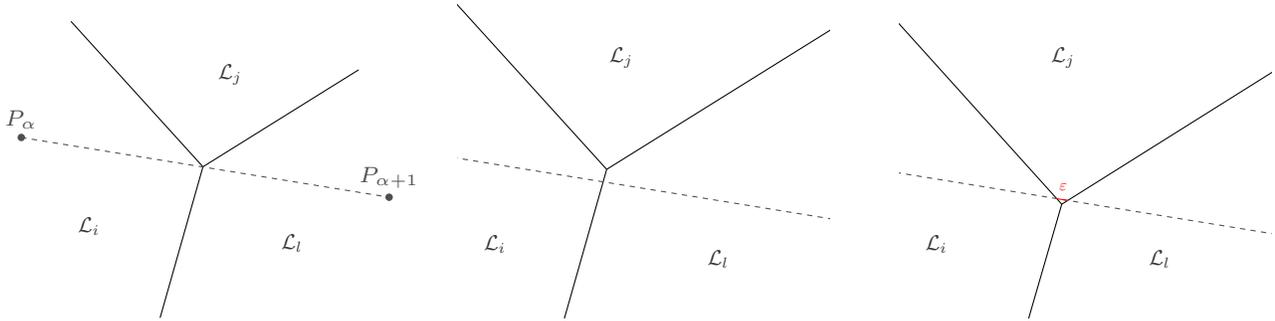
\begin{figure}[!ht]
 \begin{minipage}[t]{.32\textwidth}
 \centering
 \begin{tikzpicture}[scale=1.5]
\draw (1.54,3.07)-- (2.71,1.78);
\draw (2.71,1.78)-- (4.09,2.64);
\draw (2.71,1.78)-- (2.33,0.44);
\draw [color=black!80](1.7,1.25) node[scale=.8] {$\mathcal{L}_i$};
\draw [color=black!80](3.5,1.10) node[scale=.8] {$\mathcal{L}_l$};
\draw [color=black!80](2.95,2.6) node[scale=.8] {$\mathcal{L}_j$};
\draw[color=black!70,dash pattern=on 2pt off 2pt] (1.1,2.04)-- (4.36,1.51);
\begin{scriptsize}
\draw[fill,color=black!70] (1.1,2.04) circle   (.8pt) node[above,scale=1.2] {$P_{\alpha}$}; 
\draw[fill,color=black!70] (4.36,1.51) circle   (.8pt) node[above,scale=1.2] {$P_{\alpha+1}$}; 
\end{scriptsize}

\end{tikzpicture}
 \end{minipage}
 \begin{minipage}[t]{.32\textwidth}
 \centering
 \begin{tikzpicture}[scale=99]
\clip(2.69,1.76) rectangle (2.74,1.81);
\draw (1.54,3.07)-- (2.71,1.78);
\draw (2.71,1.78)-- (4.09,2.64);
\draw (2.71,1.78)-- (2.33,0.44);
\draw[color=black!70,dash pattern=on 2pt off 2pt] (1.1,2.04)-- (4.36,1.51);
\draw [color=black!80](2.695,1.77) node[scale=.8] {$\mathcal{L}_i$};
\draw [color=black!80](2.725,1.768) node[scale=.8] {$\mathcal{L}_l$};
\draw [color=black!80](2.712,1.795) node[scale=.8] {$\mathcal{L}_j$};
\end{tikzpicture}
 \end{minipage}
 \begin{minipage}[t]{.32\textwidth}
 \centering
 \begin{tikzpicture}[scale=99]
\clip(2.69,1.76) rectangle (2.74,1.81);
\draw (1.541663097366033,3.074111905731396)-- (2.711824949783157,1.775315239750693);
\draw (2.711824949783157,1.775315239750693)-- (4.093611818062952,2.642562853776146);
\draw (2.711824949783157,1.775315239750693)-- (2.325920509092403,0.443322492850355);
\draw[color=black!70,dash pattern=on 2pt off 2pt] (1.101730265651666,2.038160951574417)-- (4.363245216005776,1.507023656860263);
\draw[color=red] (2.711147078960937,1.776067628287411)-- (2.71263711377168,1.775824976308054);
\draw [color=red](2.712056986247088,1.776477681830871) node[scale=.6,above] {$\varepsilon$};
\draw [color=black!80](2.695,1.77) node[scale=.8] {$\mathcal{L}_i$};
\draw [color=black!80](2.725,1.768) node[scale=.8] {$\mathcal{L}_l$};
\draw [color=black!80](2.712,1.795) node[scale=.8] {$\mathcal{L}_j$};
\end{tikzpicture}
 \end{minipage}
 \caption{Example of non-uniqueness to the minimization problem \eqref{eq::nextTime}. At the left the polyline $[P_\alpha,P_{\alpha+1}]$ intersects the junction of 3 Laguerre cells $(i,j,l)$. In the center $k = l$ chosen as expected. At the right $k = j$ chosen and $t_e^{j\alpha}-t_s^{j\alpha} = \varepsilon$, a small error occurs assigning a bit of the segment to the $j$th Laguerre cell.}\label{fig::corner}
\end{figure}

There might exist several solutions $k$ to the minimization problem \eqref{eq::nextTime}. In this case, the segment encounters a corner, the intersection of at least 3 Laguerre cells.
In practice this case never occurs thanks to floating point arithmetic but the rounding error can elect a non-suitable candidate. In this case, the algorithm assigns a segment of negligible length to the candidate, see Figure~\ref{fig::corner}.

\subsection{Parallelism}

The evaluation of $g(\phi,\x)$ and its derivative with respect to $\phi$ are discussed in Section \ref{sec::integrals}. For each segment of the polyline the intersection times are computed sequentially by~\eqref{eq::nextTime} and the integration is performed on the fly. The parallelism is unfurled at the highest level of the algorithm for the segments. Remember that the computation of the intersection of the polyline and the Laguerre tessellation requires for each segment the knowledge of the cell of its starting point. This cell can be inferred from the calculation of the previous segment. Hence, we slice the polyline into contiguous chunks of segments with equivalent size, each worker dealing with one chunk.
For each of these chunks the sequential integration procedure is launched and then the results are merged. We carry out both OpenMP and \texttt{c++11} thread class implementation. For a standard chunk size and for a regular processor the \texttt{c++11} performance overtakes the OpenMP's one. We suppose that this difference of performance (up to a factor two) is due to OpenMP overhead.

\medskip

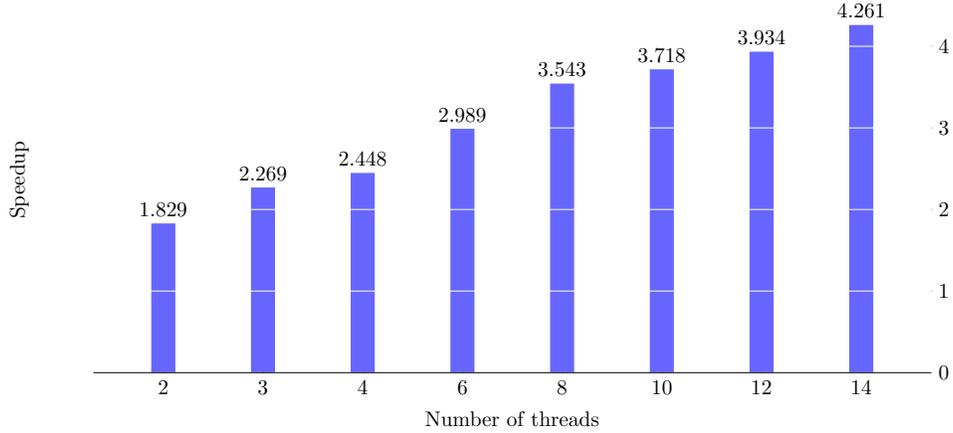
\begin{figure}
 \centering
 \begin{tikzpicture}[scale=0.8]
  \centering
  \begin{axis}[
        ybar, axis on top,
        height=8cm, width=15.5cm,
        bar width=0.4cm,
        ymajorgrids, tick align=inside,
        major grid style={draw=white},
        enlarge y limits={value=.1,upper},
        ymin=0, ymax=4.3,
        axis x line*=bottom,
        axis y line*=right,
        y axis line style={opacity=0},
        tickwidth=0pt,
        enlarge x limits=true,
        legend style={
            at={(0.5,-0.2)},
            anchor=north,
            legend columns=-1,
            /tikz/every even column/.append style={column sep=0.5cm}
        },
        xlabel={Number of threads},
        ylabel={Speedup},
        symbolic x coords={
           2,3,4,6,8,10,12,14},
       xtick=data,
       nodes near coords={
        \pgfmathprintnumber[precision=3]{\pgfplotspointmeta}
       }
    ]
    \addplot [draw=none, fill=blue!60] coordinates {    							
      (2,1.8290033222)
      (3, 2.269405933)
      (4,2.4480260321)
      (6,2.9888826614) 
      (8,3.543337456) 
      (10,3.7177647102)
      (12,3.9339883116)
      (14,4.2609571393)};
  \end{axis}
\end{tikzpicture}
 \caption{Speedup for the computation of the cost function $g$ with \texttt{c++11} thread class. The parallelism is placed over the segment integration for $n=200\mathrm{K}$ points and $p=80\mathrm{K}$ segments for an increasing number of threads with a super-calculator equipped with Intel Xeon\textsuperscript{\textregistered} E5-2680}\label{fig::speedup}
\end{figure}

In Figure~\ref{fig::speedup}, the evolution of the computation time of the cost function $g$ and its gradient with respect to $\phi$ is benchmarked. The speedup unit is defined as the ratio of the execution time of the non-parallelized task over the parallelized task. OpenMP and \texttt{c++11} thread class are both implementations of shared memory parallelization, the memory is simultaneously accessible for every thread. The performance of a process depends intrinsically on how close the data is. If the data is on processor cache the latency is a dozens of CPU cycles, if the data is located on RAM, the latency is $40ns$. In our experiment, for large number of threads, $30$\% of the data was not located on cache. This explains why performance drops as the number of threads increases.

\section{Computation of the optimal transport}
\label{sec::four}

\begin{figure}[!ht]
 \begin{minipage}[t]{.5\textwidth}
 \centering
 \begin{tikzpicture}[scale=1.1]
\definecolor{xdxdff}{rgb}{0.49,0.49,1}
\definecolor{uuuuuu}{rgb}{0.27,0.27,0.27}
\definecolor{qqqqff}{rgb}{0,0,1}
\draw (3,0.8)-- (3.01,-1.2);
\draw (2.99,3.35)-- (2.99,5.35);
\draw (0,3.79)-- (1.89,2.71);
\draw (1.9,1.43)-- (0,0.33);
\draw (4.1,1.44)-- (6,0.35);
\draw (4.1,2.71)-- (6,3.82);
\draw (1.9,1.43)-- (3,0.8);
\draw (4.1,1.44)-- (3,0.8);
\draw (1.89,2.71)-- (2.99,3.35);
\draw (2.99,3.35)-- (4.1,2.71);
\draw (1.89,2.71)-- (1.9,1.43);
\draw (4.1,2.71)-- (4.1,1.44);
\draw[color=black!70,dash pattern=on 2pt off 2pt] (1.91,4.04)-- (4.36,3.53);
\draw[color=black!70,dash pattern=on 2pt off 2pt] (4.36,3.53)-- (5.64,2.04);
\draw[color=black!70,dash pattern=on 2pt off 2pt] (5.64,2.04)-- (4.48,0.35);
\draw[color=black!70,dash pattern=on 2pt off 2pt] (4.48,0.35)-- (0.43,1.02);
\begin{scriptsize}
\draw[fill] (3,2.07) circle   (.8pt) node[below,scale=1.2] {$x_i$};
\draw[fill,color=black!70] (1.91,4.04) circle   (.8pt) node[above,scale=1.2] {$P_{\alpha-2}$};
\draw[fill,color=black!70] (4.36,3.53) circle   (.8pt) node[above,scale=1.2] {$P_{\alpha-1}$};
\draw[fill,color=black!70] (0.43,1.02) circle   (.8pt) node[left,scale=1.2] {$P_{\alpha+2}$};
\draw[fill,color=black!70] (4.48,0.35) circle   (.8pt) node[below,scale=1.2] {$P_{\alpha+1}$};
\draw[fill,color=black!70] (5.64,2.04) circle   (.8pt) node[right,scale=1.2] {$P_{\alpha}$};
\end{scriptsize}
\end{tikzpicture}
 \end{minipage}
 \begin{minipage}[t]{.49\textwidth}
 \centering
 \pgfplotsset{cellmodel/.style={%
width=1.1\linewidth,
axis lines=center,
axis equal image,
domain=0:800,
xmin= 0, xmax= 7,
ymin=-1.5, ymax=7,
restrict x to domain=0:6.2,
ticks=none,
}}

\begin{tikzpicture}[scale=1]
\begin{axis}[cellmodel,
        xlabel=$\phi_i$,
        ylabel=$\nabla g(\phi_i)$]
        
    \addplot[black] plot coordinates {
        (0,6)
        (2,6)
    };
    \addplot[mark=*,black] plot coordinates {
        (2,6)
        (3,5.5)
        (3.8,4.7)
        (4.8,2.7)
    };
    \addplot[black] plot coordinates {
        (4.8,2.7)
        (6.1,-1.3)
    };
\end{axis}
\end{tikzpicture}
 \end{minipage}
 \caption{Counter-example to the smoothness of $\nabla g$. Consider the a Laguerre tessellation and a set of segment as displayed in the left and increase $\phi_i$. The Laguerre cell $\lag_i$ increases and the gradient of $g$ with respect to $\phi$ exhibits kinks. Those kinks happen each time a Laguerre cell meets a segment for the first time. Note that under Hypothesis \eqref{eq::generic}, the boundary of the Laguerre cell can only cross transversally the segment. If Hypothesis \eqref{eq::generic} is not met $\nabla g$ may fail to exist.}\label{fig::notC2}
\end{figure}
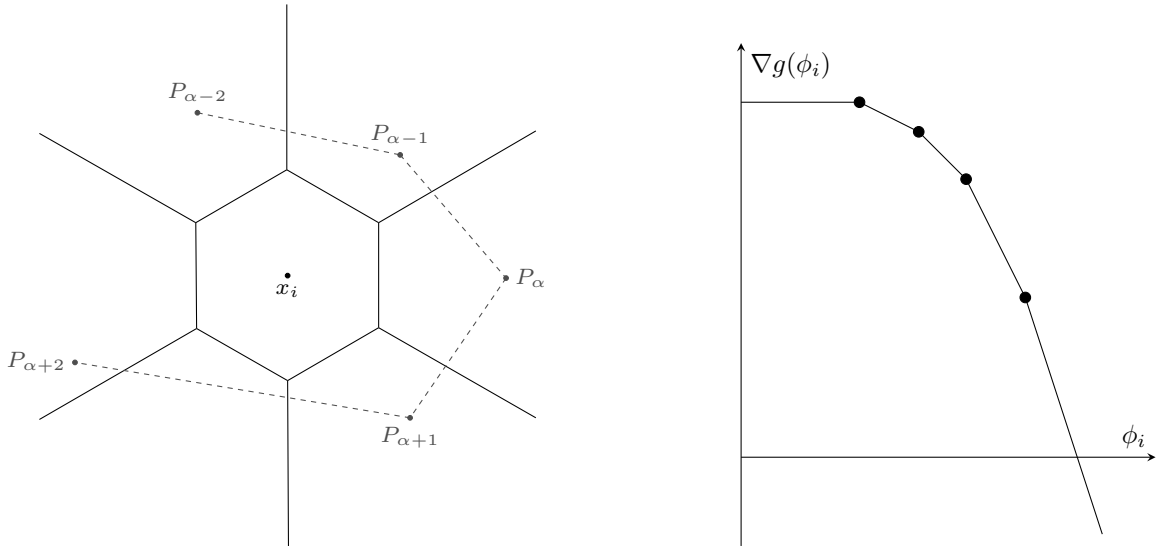

The goal of this section is to compare the available methods to optimize $g$ with respect to $\phi$ \eqref{eq::costfunc}. The functional $g$ is concave with respect to $\phi$ since it is a dual formulation of the problem \eqref{eq::MongeKanto} see \cite{villani2008optimal}. Under the hypothesis \eqref{eq::generic} the functional $g$ is $\mathcal{C}^1$ with respect to $\phi$. Note that a higher level of regularity cannot be established as argued in Figure \ref{fig::notC2}. This lack of regularity precludes the use of the convergence framework for second order methods established by M\'erigot \textit{et al.} \cite{kitagawa2016convergence}.

However, we may apply classical convergence results of concave analysis to functional the $g$. To trigger those arguments, the Lipschitz constant of the gradient $L$ has to be bounded. Under Hypothesis \eqref{eq::generic}, using the notation of Section~\ref{sec::Hess} summarized in Figure \ref{fig::displacementBoundary} and the expression of the Hessian matrix \eqref{eq::Hessian}, we establish a finite upper bound for $L$ using the Gershgorin circle theorem : 

\begin{equation} \label{eq::LipschitzConstant}
 | L |  \leq 2 \max_i \sum_{j\neq i} \left| \frac{\partial^2 g }{\partial \phi_i \partial \phi_j} \right| \leq  \frac{p^2 \quad \displaystyle{\max_\alpha} (\rho_\alpha)}{\displaystyle{\min_{l \neq m}} (\| \x_l - \x_m \|_2) \quad  \displaystyle{\min_{i,j,\alpha,i\neq j}} | \cos((\theta^{ij}_\alpha)) | },
\end{equation}
where $p$ is the number of segments composing the polyline, $\rho_\alpha$ is the density associated to each segment and $\theta_\alpha^{ij}$ is the angle between the vector $P_{\alpha+1} - P_\alpha$ and $(\x_j-\x_i)$. For fixed measures $\nu$ and $\mu$ the evaluation of the bound~\eqref{eq::LipschitzConstant} is costly as it involves combinatorial quantities and is not evaluated in practice. 

The gradient Lipschitz condition \eqref{eq::LipschitzConstant} is sufficient to ensure convergence of ascent methods. Note that $g$ is bounded from above as $\mu$ and $\nu$ are compactly supported measures. The gradient method with step $s < \frac 2 L$ converges to a stationary point $\phi^\star$ see \cite{polyak1987introduction,bertsekas1999nonlinear,nesterov2013introductory}. The same holds for variant step-size method with line search such as Global Brazalai Borwein algorithm \cite{raydan1997barzilai,fletcher2005barzilai}.

Because of the counter-example in Figure~\ref{fig::notC2}, and the lack of regularity of the gradient of $g$, quadratic convergence cannot be guaranteed for second order methods. Note also that in a generic setting, the Newton method is impractical since the Hessian matrix defined in equation~\eqref{eq::Hessian} fails to be invertible, see Section~\ref{sec::HessMethod}.

\subsection{Choice of optimization method}
The different methods are tested against the same benchmark. It consists in drawing uniformly $10\mathrm{K}$ points and $500$ segments in 2D. Several realization of the optimization procedure are performed and plotted in the corresponding figures. The methods are benchmarked in Table \ref{tab:comparison}.
\subsubsection{First order method}

\begin{figure}[!ht]
 \begin{minipage}[t]{.5\textwidth}
 \centering
 \includegraphics[width=\textwidth]{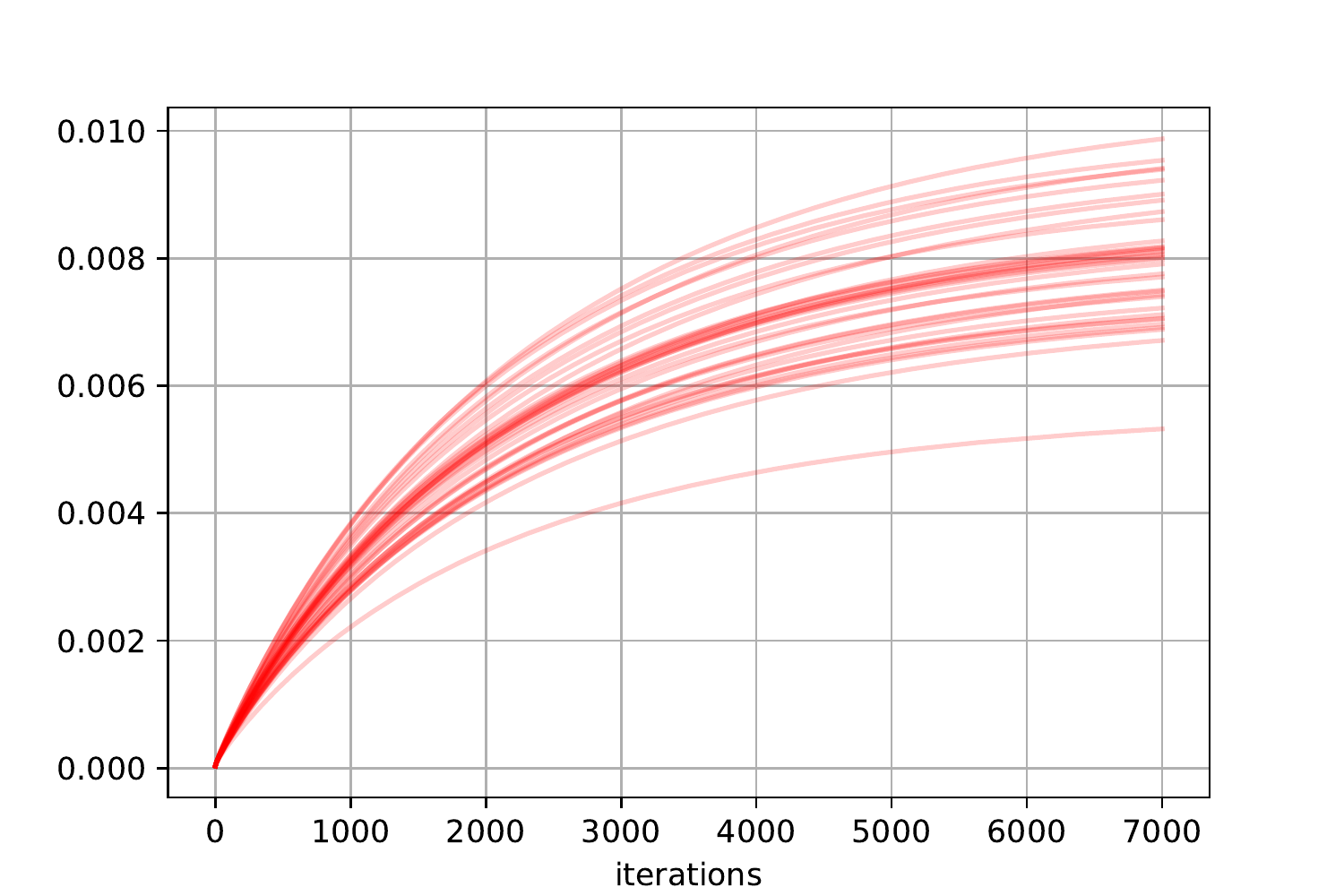}
 \end{minipage}
 \begin{minipage}[t]{.49\textwidth}
 \centering
 \includegraphics[width=\textwidth]{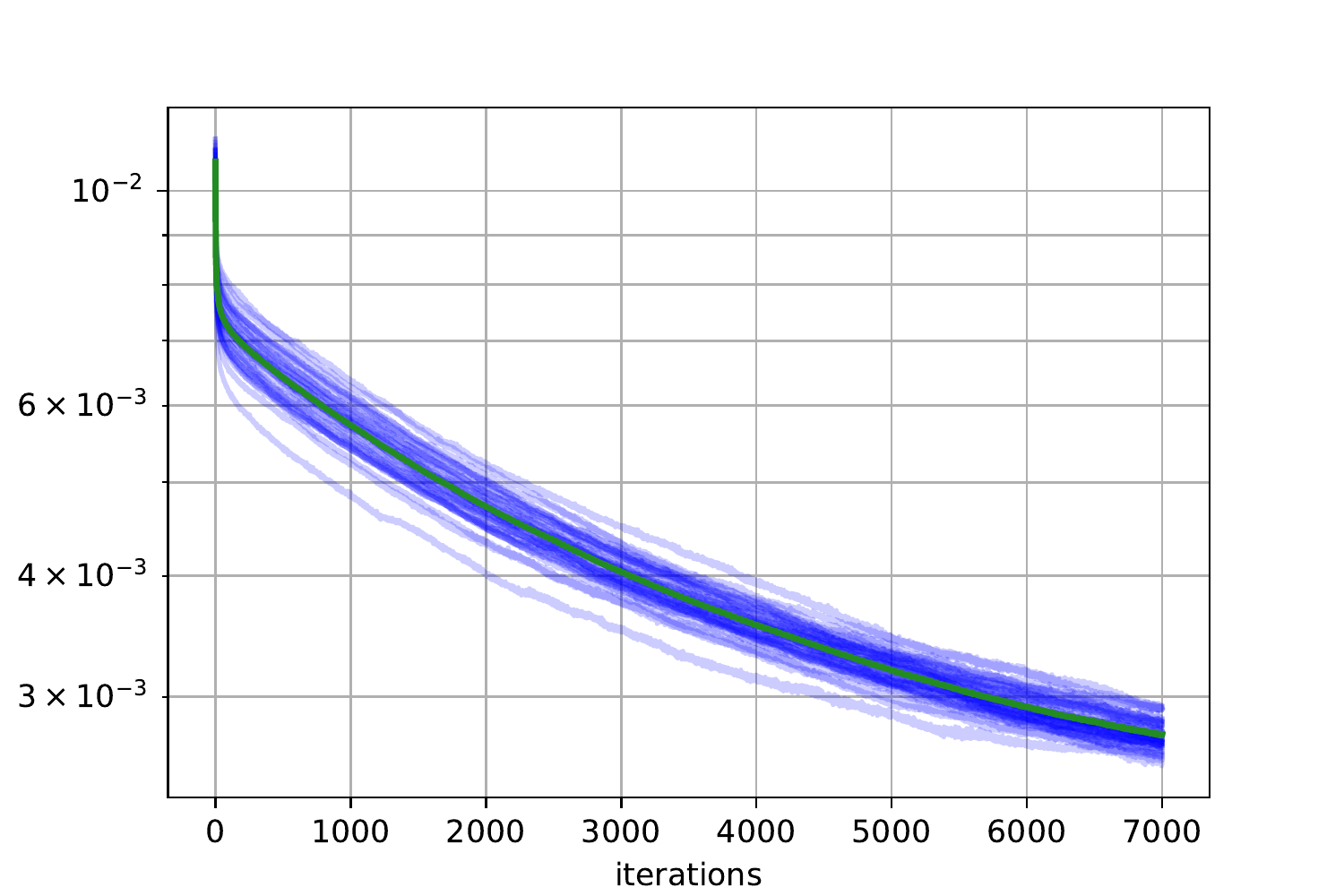}
 \end{minipage}
 \caption{Gradient ascent method : cost function history (left) and norm of the gradient (right). Method stopped by the max iteration criterion.}\label{fig::convergenceAscnt}
\end{figure}

We first implement a gradient ascent method. The generic convergence history is displayed in Figure~\ref{fig::convergenceAscnt}. As explained in Section~\ref{sec::four}, the Lipschitz constant of the gradient cannot be satisfactorily computed, so the step-size is chosen according to strong Wolfe conditions. In practice the algorithm settles for a constant step of $0.05$. As it can be observed in Figure~\ref{fig::convergenceAscnt} the rate of convergence of the gradient towards zero is too slow for this method to be used in practice.

\begin{figure}[!ht]
 \begin{minipage}[t]{.5\textwidth}
 \centering
 \includegraphics[width=\textwidth]{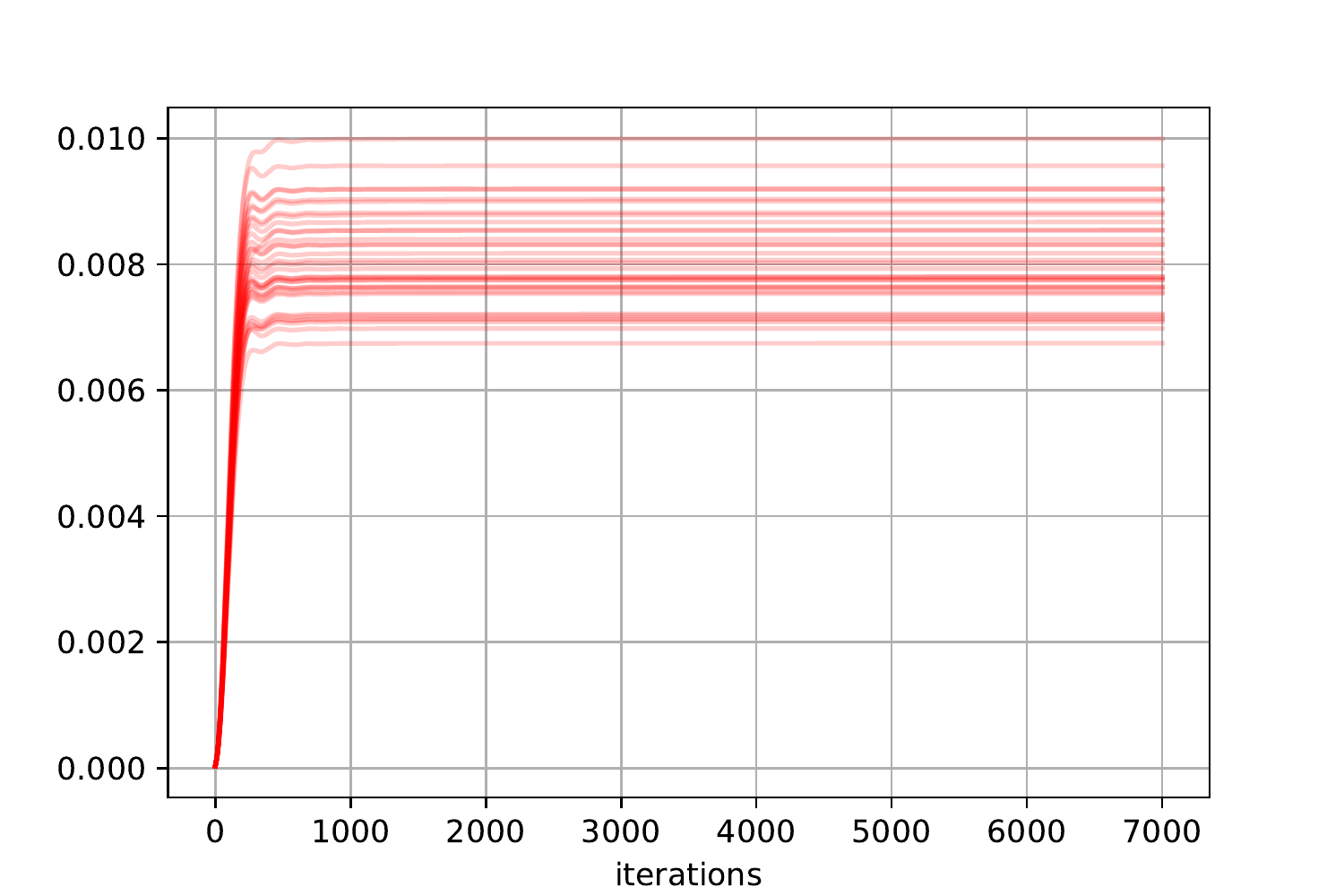}
 \end{minipage}
 \begin{minipage}[t]{.49\textwidth}
 \centering
 \includegraphics[width=\textwidth]{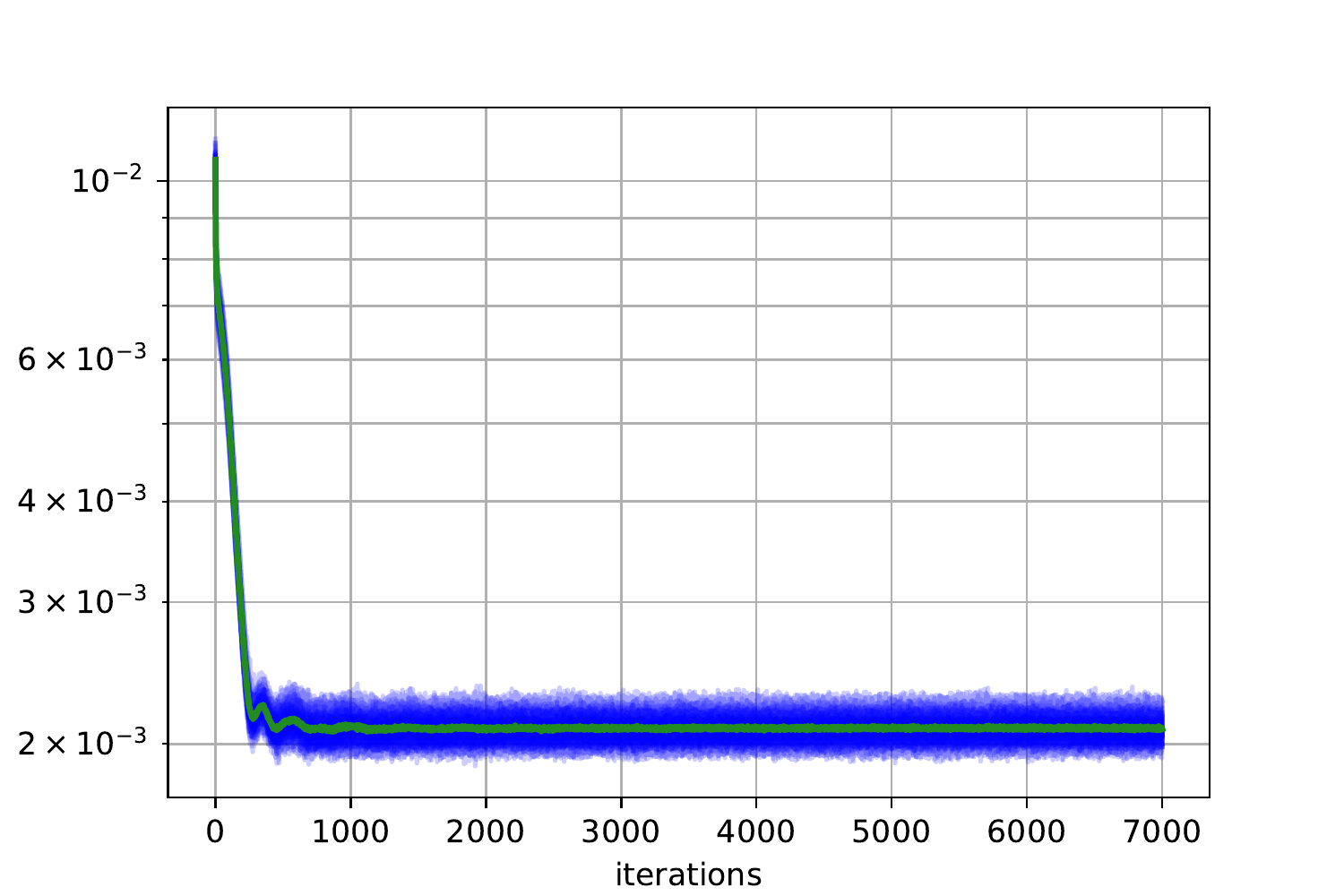}
 \end{minipage}
 \caption{Nesterov acceleration of gradient ascent :  cost function history (left) and norm of the gradient (right). Method stopped by the max iteration criterion.}\label{fig::convergenceNeste}
\end{figure}

We have tried other first order methods including Polak-Ribière, Fletcher-Reeves, Barzalai Borwein and Nesterov acceleration. Only Nesterov acceleration has behaved differently from the other methods and it is displayed in Figure~\ref{fig::convergenceNeste}. Note however, that Nesterov acceleration requires an estimation of $L$, the Lipschitz constant of the gradient. Figure~\ref{fig::convergenceNeste} was obtained with an estimation of $L$ that promotes quick decreases of gradient norm at the beginning of the algorithm. Note however that the estimation of $L$ seems too optimistic since the Nesterov algorithm does not converge. A more pessimistic choice of $L$ leads to a convergence rate similar to the gradient algorithm, at least during the first $7\mathrm{K}$ iterations.

\subsubsection{Quasi-Newton Method}

\begin{figure}[!ht]
 \begin{minipage}[t]{.5\textwidth}
 \centering
 \includegraphics[width=\textwidth]{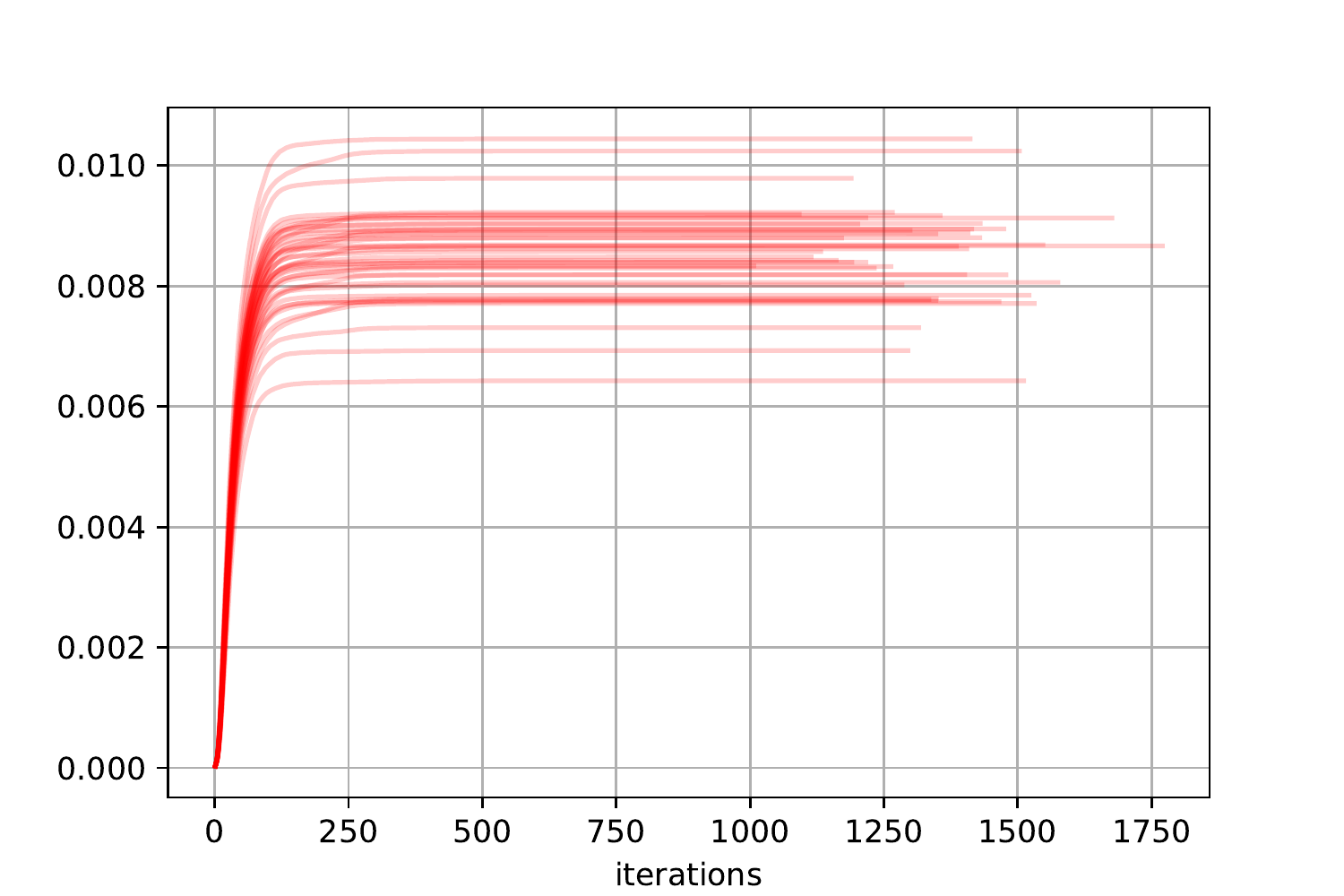}
 \end{minipage}
 \begin{minipage}[t]{.49\textwidth}
 \centering
 \includegraphics[width=\textwidth]{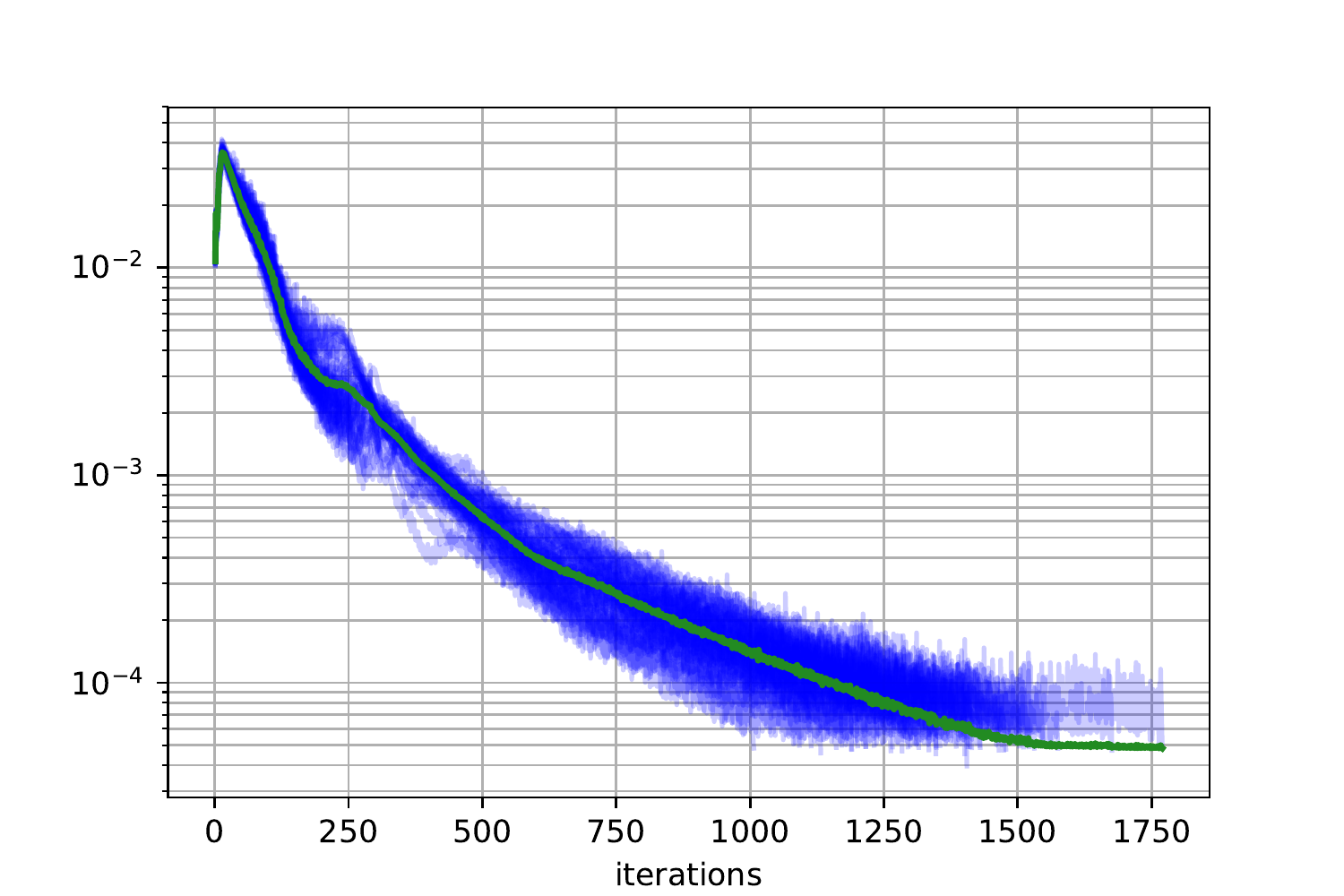}
 \end{minipage}
 \caption{L-BFGS method : cost function history (left) and norm of the gradient (right). Method converges to the targeted gradient norm $5.10^{-5}$ in a average of 1500 iterations.}\label{fig::convergenceBFGS}
\end{figure}

The main idea behind quasi-Newton methods is to build at the iteration $k$ an approximation of the Hessian matrix of $g$. We choose the limited memory BFGS method, which only stores a limited amount of vectors determined by the user. The result is displayed in Figure~\ref{fig::convergenceBFGS}. This method converges faster than first order method, the L-BFGS algorithm reaches the desired gradient tolerance $5.10^{-5}$ within $1500$ iterations. By contrast first order methods failed to converge in 7000 iterations. However when the method is close to a critical point the convergence speed tends to decrease.

\medskip

\subsubsection{Newton and quasi-Newton hybridization}
\label{sec::HessMethod}
In the light of the previous section, the flaw of the quasi-Newton method is its lack of briskness when it reaches the  vicinity of a maximizer. Second-order methods are known to converge quadratically in the basin of attraction of a maximizer. The idea of the hybrid algorithm is to start with a quasi-Newton method, and then switch to a second order Newton algorithm when the basin of attraction is reached. The quandary is to determine when to start the Newton method in preference to the L-BFGS method. Indeed far from the basin of attraction the direction given by the Newton method is poor and the natural step $1$ is truncated by the line search which causes extra computation time.

\begin{figure}[!ht]
 \begin{minipage}[t]{.5\textwidth}
 \centering
 \includegraphics[width=\textwidth]{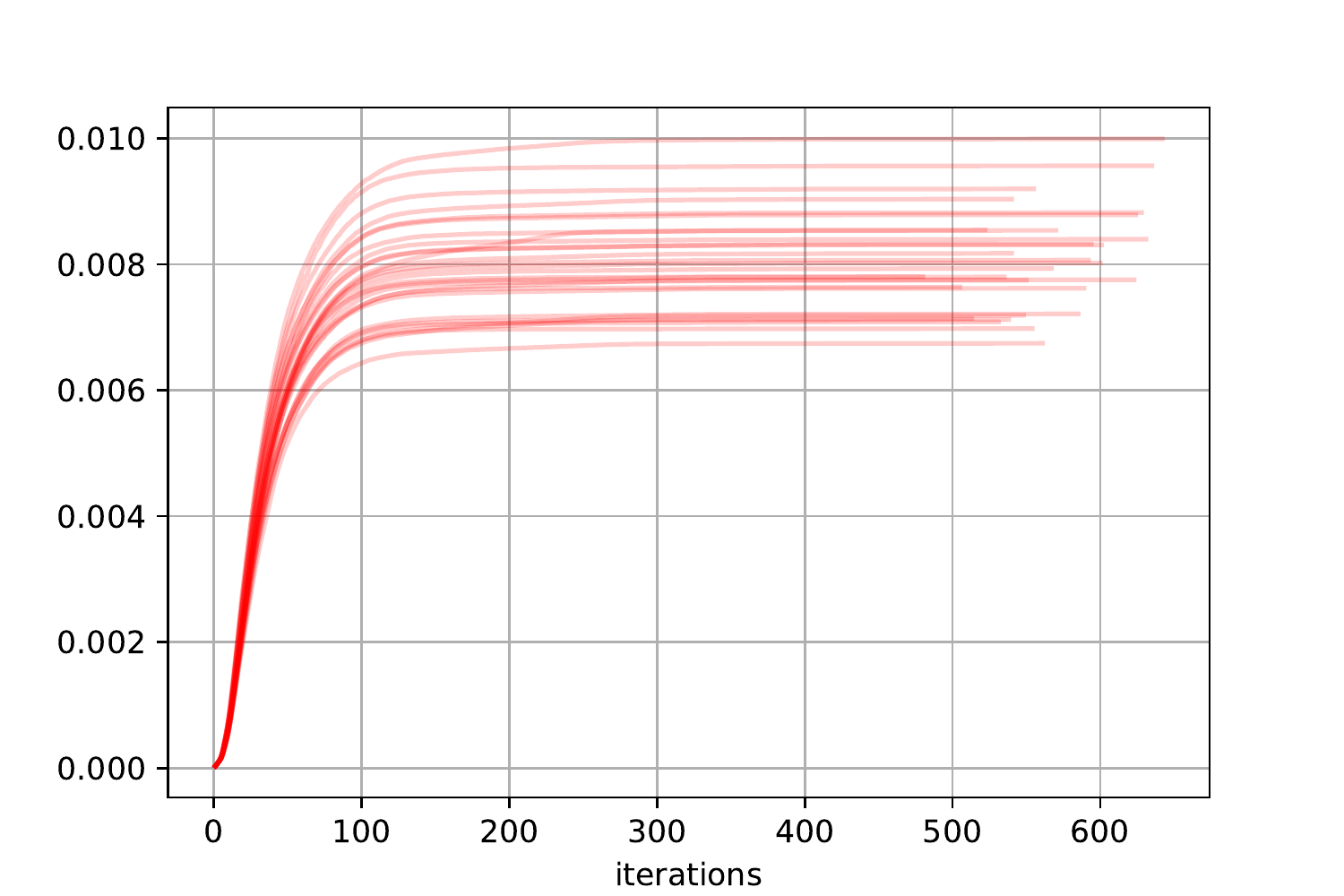}
 \end{minipage}
 \begin{minipage}[t]{.49\textwidth}
 \centering
 \includegraphics[width=\textwidth]{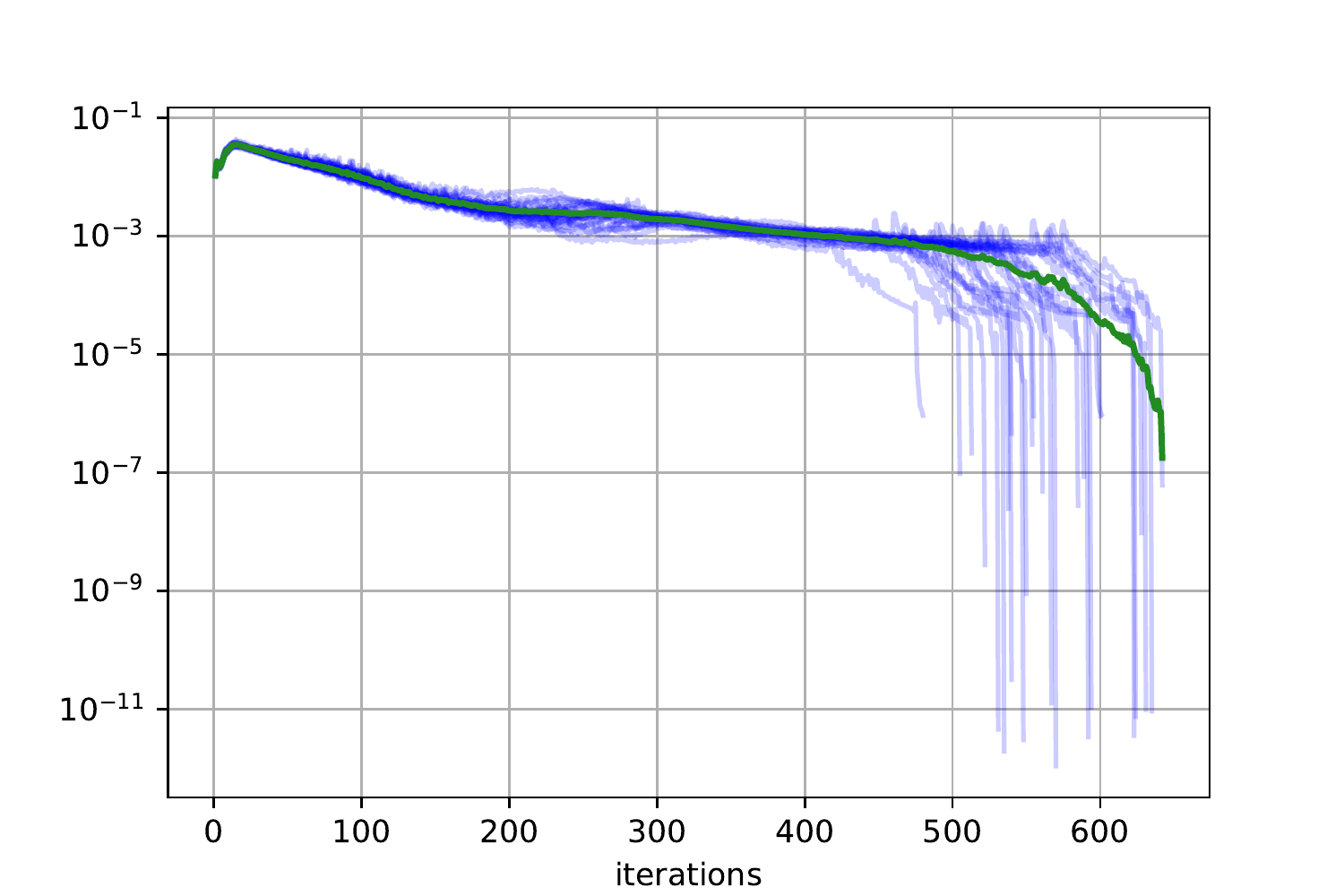}
 \end{minipage}
 \caption{Hybrid Newton method : cost function history (left) and norm of the gradient (right). The method converges to the targeted gradient norm $10^{-6}$ in 600 iterations on average.}\label{fig::NewtonMethod}
\end{figure}

We propose to switch from the L-BFGS algorithm to the Newton algorithm as soon as there is no empty Laguerre cell. A Laguerre cell $\lag_i$ is called empty as long as there is no mass assigned to its centroid $\x_i$ that is $\nu(\lag_i) = 0$. This criterion is informally motivated by the fact that the second order information relates the competition and the connectivity between Laguerre cells, hence if one of them is empty its information is inconsequential. From a practical point of view, this condition ensures the Hessian on being a full rank matrix and hence invertible. This procedure is described in Algorithm \ref{alg::GradientPhiDescent} and the result are displayed in Figure~\ref{fig::NewtonMethod}.

\begin{algorithm}
\caption{Computation of optimal transport}\label{alg::GradientPhiDescent}
\begin{algorithmic}[1]
\Ensure Dirac positions $\x$ and masses $m$.
\Ensure Nodes $P$ and density $\rho$ of the polyline.
\Require $\phi_{\text{init}}$ a starting $\phi$ for the computation of $g$
\Require gradTol L2 tolerance on the gradient norm
\Require outerMax maximum number of iterations 
\Require $\texttt{wolfeProcedure}$ standard line search with strong Wolfe condition and initial step $s=1$. \cite{bertsekas1999nonlinear}
\Function{$\texttt{computeOptimalTransport}$}{$\x$,$m$,$P$,$\rho$}
\State $\texttt{bf} \gets \texttt{LBFGS}(\text{memSize})$ \Comment{Initialization of L-BFGS}
\State $ \phi \gets \phi_{\text{init}}$
\State $\nabla \phi$,cost,hiddenNumber $\gets \texttt{computeIntegration}(\phi)$
\State $i \gets 0$
\While {$i < \text{outerMax} \And \text{gradTol} < \|\nabla \phi \|_2$}
  \If {hiddenNumber $\neq 0$} 
    \State d $\gets \texttt{bf.findDirection}(\nabla \phi)$ 
  \Else
    \State $\Hess \gets \texttt{computeHessian}(\phi)$
    \State d $\gets - \Hess^{-1} \nabla \phi$ \Comment{The Hessian is definite so Newton direction is taken}
  \EndIf
  \State $s \gets \texttt{wolfeProcedure}(\text{cost},\phi,s,\text{d})$
  \State $\nabla \phi_{\text{old}} \gets \nabla \phi$
  \State $\nabla \phi$,cost,hiddenNumber $\gets \texttt{computeIntegration}(\phi+s \text{d})$
  \State $\texttt{bf.addDirection}(\nabla \phi_{\text{old}} - \nabla \phi,sd)$ \Comment{Actualize memory of the L-BFGS algorithm}
  \State $\phi \gets \phi + s d$
\EndWhile
\EndFunction

\end{algorithmic}
\end{algorithm}

In the numerical tests, the switch between the two Newton method happens approximatively $30$ iterations before termination of the algorithm. During the Newton phase of the hybrid algorithm, the algorithm undergoes two stages. During the twenty first iteration of the Newton method the algorithm stabilizes around the maximum and the Wolfe line-search prevents picking $s=1$. During the second stage, quadratic convergence is achieved in few iterations.

Levenberg-Marquartdt's method was also implemented with a regularization parameter that tends to zero when approaching the basin of convergence, however in large scale optimization problems the lightweight hybrid method is preferable since it requires the inversion of the Hessian only in the last iterations. 

\begin{table}[H]
   \small 
   \centering 
   \begin{tabular}{l|ccccc} 
   \toprule[\heavyrulewidth]\toprule[\heavyrulewidth]
   \textbf{} & \textbf{BBG} & \textbf{Nesterov} & \textbf{BFGS} & \textbf{BFGS/Newt} & \textbf{LM}\\ 
   \midrule
   Time/iteration & $36.4 m s$ & $36.7 m s$ & $49.4 m s$ & $186 m s$& $649 m s$\\
   gradient norm & $5.5 \  10^{-3}$ & $1.2 \ 10^{-3}$ & $1.5 \ 10^{-4}$ & $10^{-15}$& $5.4 \ 10^{-5}$\\
   \bottomrule[\heavyrulewidth] 
   \end{tabular}
   \caption{Comparison of algorithms for solving the optimal transport problem for $10^4$ points and $500$ lines. Time required for an iteration in millisecond and the gradient norm after 1000 iterations. The BFGS/Newton method converges up to numerical error.}\label{tab:comparison}
\end{table}

\section{Numerical examples}
\label{sec::five}
\subsection{Representation of picture}
\begin{figure}[!ht]
 \includegraphics[width=0.50\textwidth]{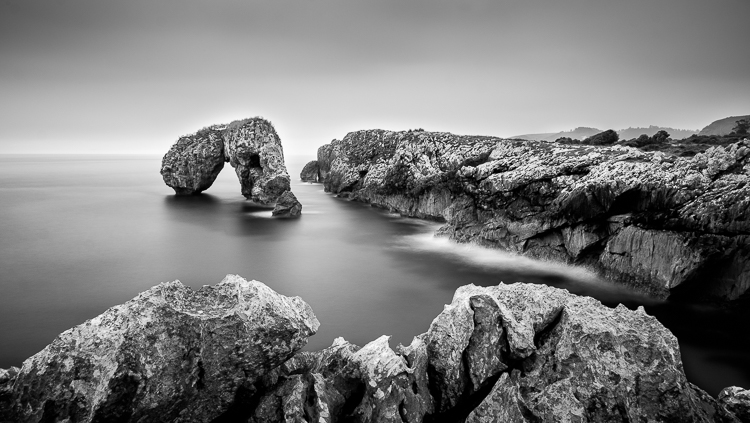}
 \includegraphics[width=0.50\textwidth]{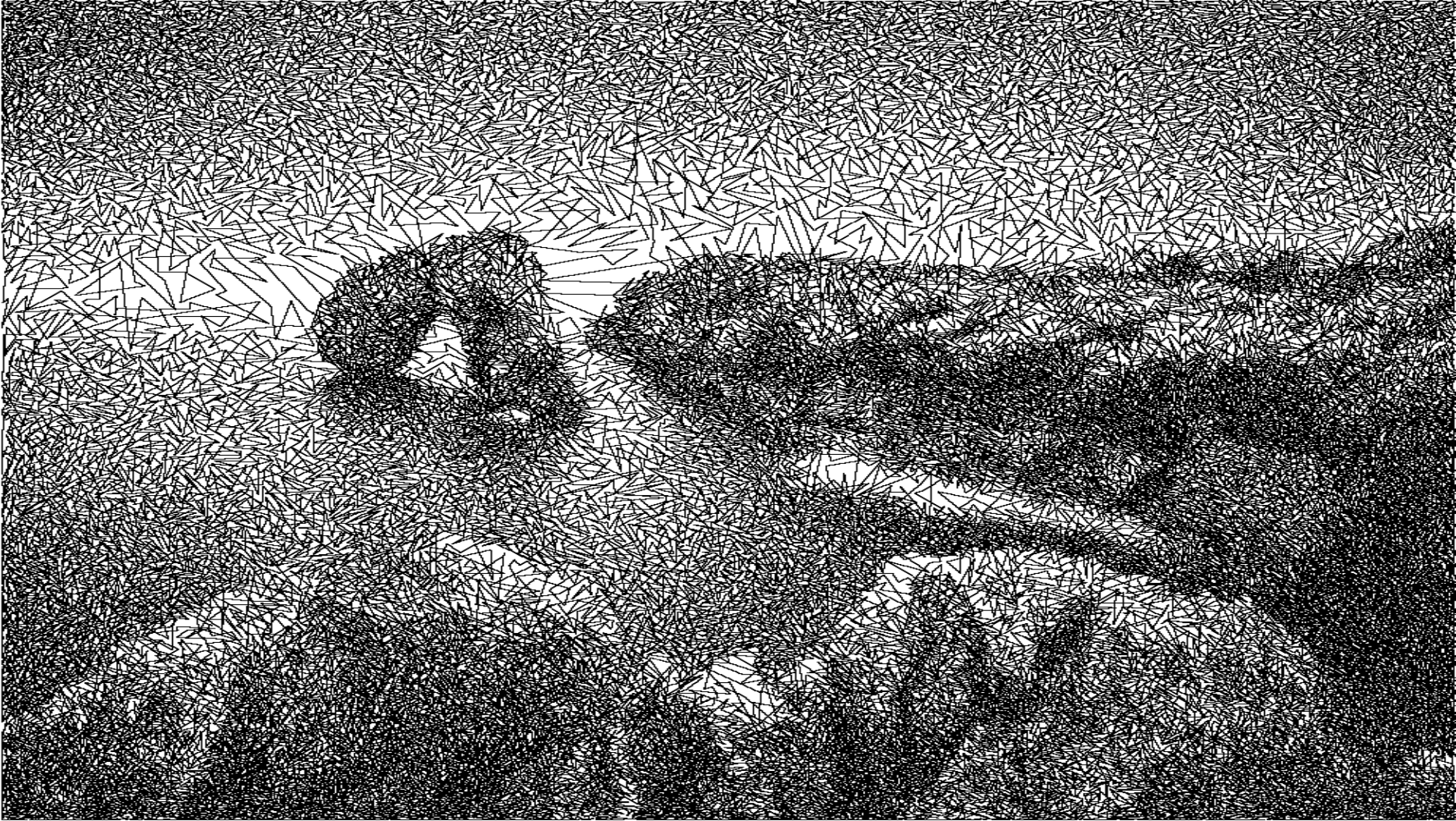}
 \includegraphics[width=0.50\textwidth]{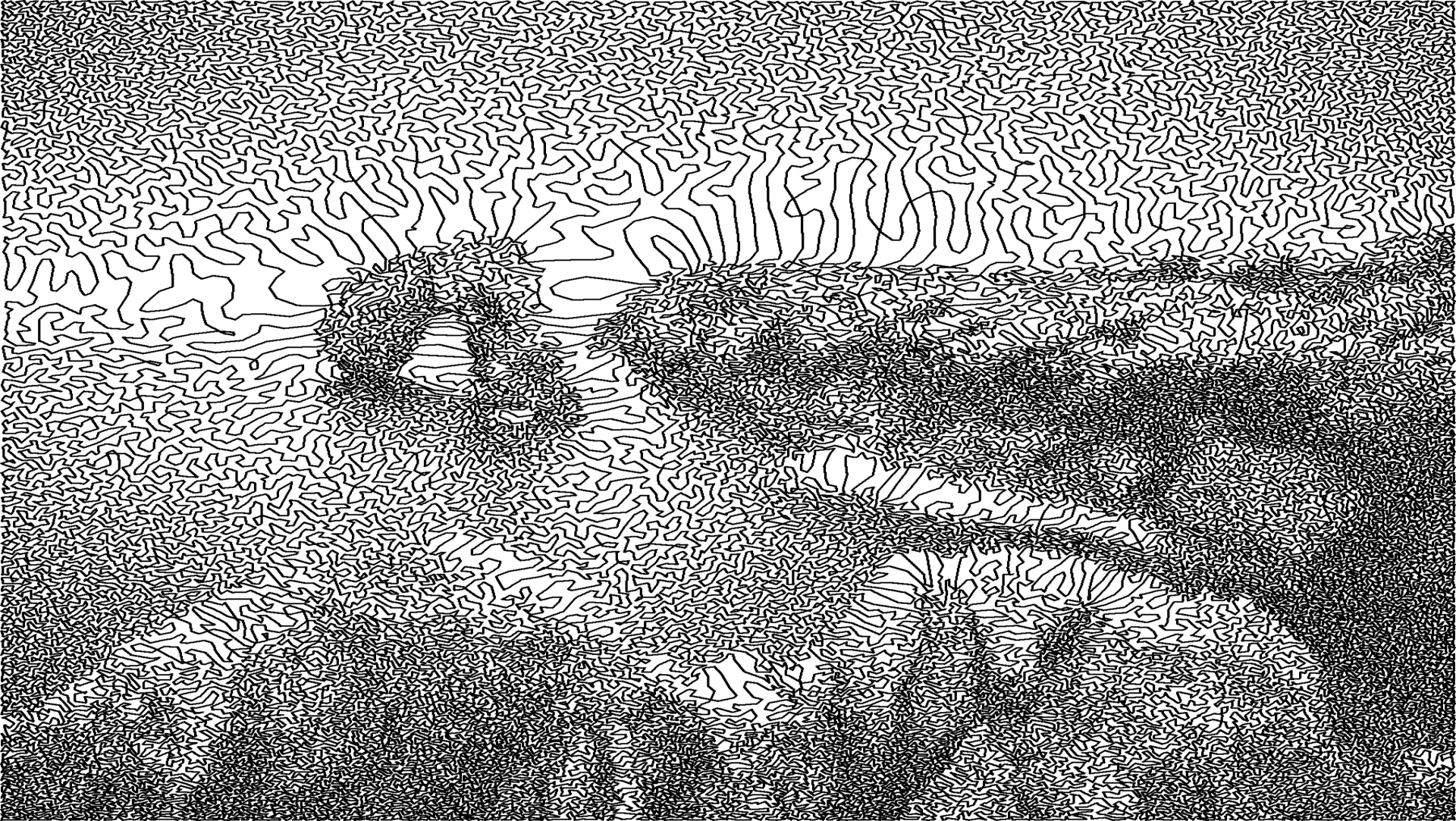}
 \includegraphics[width=0.50\textwidth]{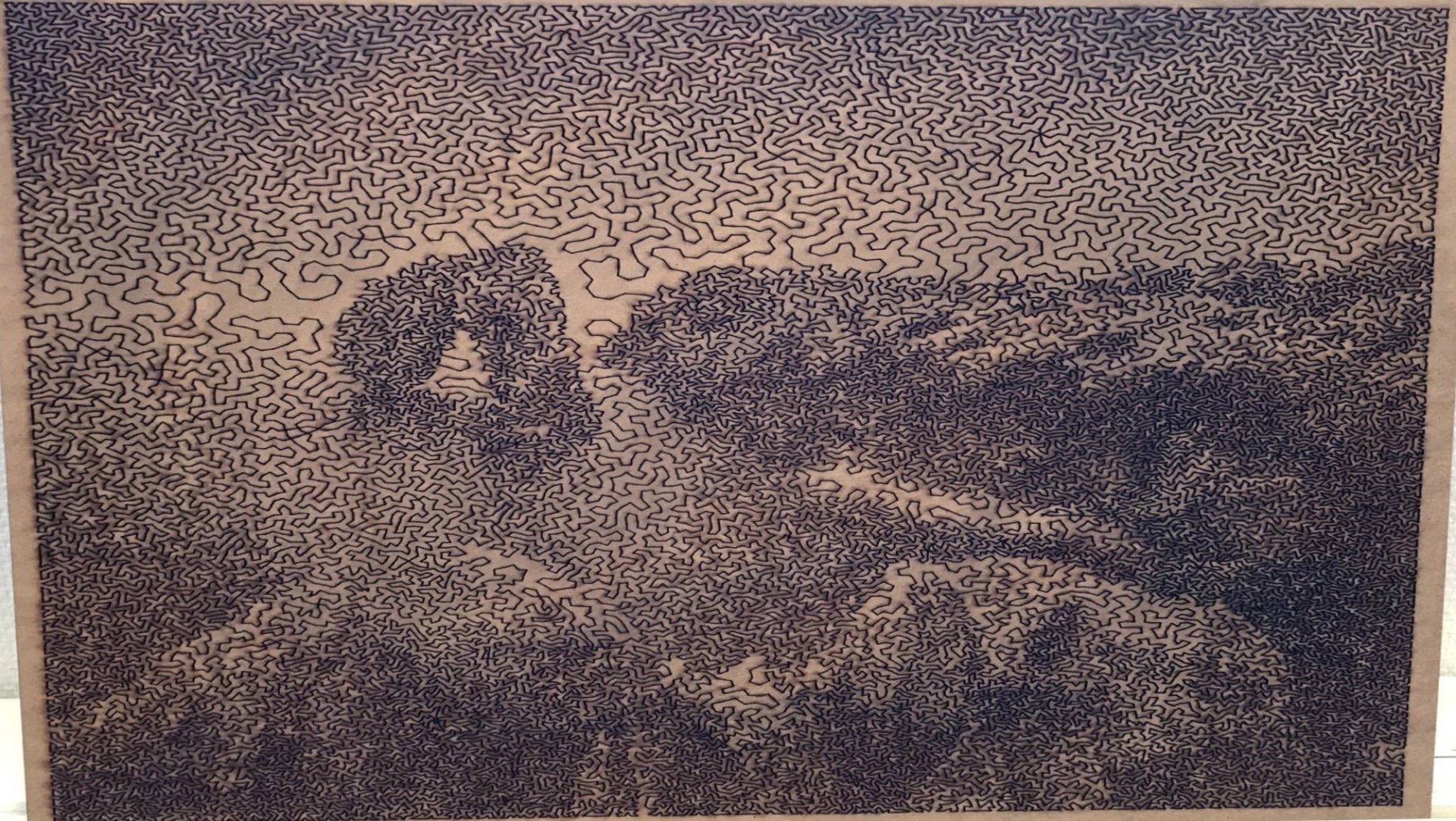}
 \caption{Curvling by $\frac 3 4$-discrete optimal transport. Original image (top left), approximation by a polyline with $400\mathrm{K}$ segments without any constraints (top right), approximation by a polyline with kinematic constraints and $250\mathrm{K}$ segments (bottom left), final rendering after wood engraving, the polyline is composed of $80\mathrm{K}$ segments (bottom right). In all tests  measure $\mu$ is represented by $320\mathrm{K}$ Dirac masses.}\label{fig::Result2D}
\end{figure}

The approximation of a measure by a curve has multiple applications of which some are described in \cite{Fred2018Projections}. In this section we discuss the representation of a picture by a polyline. The first step is the discretization of the picture, here a landscape, by a sum of Dirac masses of weights $m$. Several approaches are possible but the most intuitive is to take the Dirac positions $\x$ on a Cartesian grid and their weights equal to the pixels intensities. The polyline is then initialized randomly, and Algorithm~\ref{alg::GradientDescent} is launched. With no speed or curvature constraints, the gradient method in $P$ described in Algorithm~\ref{alg::GradientDescent} empirically gives a stationary point. This solution is highly non-smooth, the length of the segments and the angle between them are disparate see Figure~\ref{fig::Result2D} top right.

This scribbled solution is, in practice, difficult to carry out with laser engravers. Indeed for the same number of segments $p$, the engraving time of uncontrolled trajectories can be up to five times longer. In order to get around this problem, we follow the method described in \cite{Fred2018Projections}, that is projecting the polyline after the gradient step of Algorithm~\ref{alg::GradientDescent} on a set of kinematic constraints. The constraints imposed on the speed and the acceleration polyline read as : 
\begin{equation}\label{eq::kinematicConstraints}
 \mathcal{K}(K_1,K_2) = \left\lbrace (P_\alpha)_{\alpha \in \llbracket 1,p \rrbracket} \text{ s.t } \left| \begin{array}{ll}
		     \| P_{\alpha+1} - P_{\alpha} \| \leq K_1, & \forall \alpha \in \llbracket 1, p-1 \rrbracket \\               
                     \| 2P_{\alpha} - P_{\alpha-1}- P_{\alpha+1} \| \leq K_2, & \forall \alpha \in \llbracket 2, p-1 \rrbracket
                                                 \end{array}\right.\right\rbrace,
\end{equation}
with $K_1$ the constant controlling the speed of the polyline and $K_2$ its acceleration. The projection on the set $\mathcal{K}(K_1,K_2)$ is performed using the Alternating Direction Method of Multipliers (ADMM) \cite{glowinski2014alternating}.

The authors wish to thank Andrew Gibson for letting them use the original picture and Alban Gossard for his priceless help in realizing the wood engraving.

\subsubsection{Galaxy filaments}

\begin{figure}[!ht]
 \centering
 \includegraphics[width=.32\textwidth]{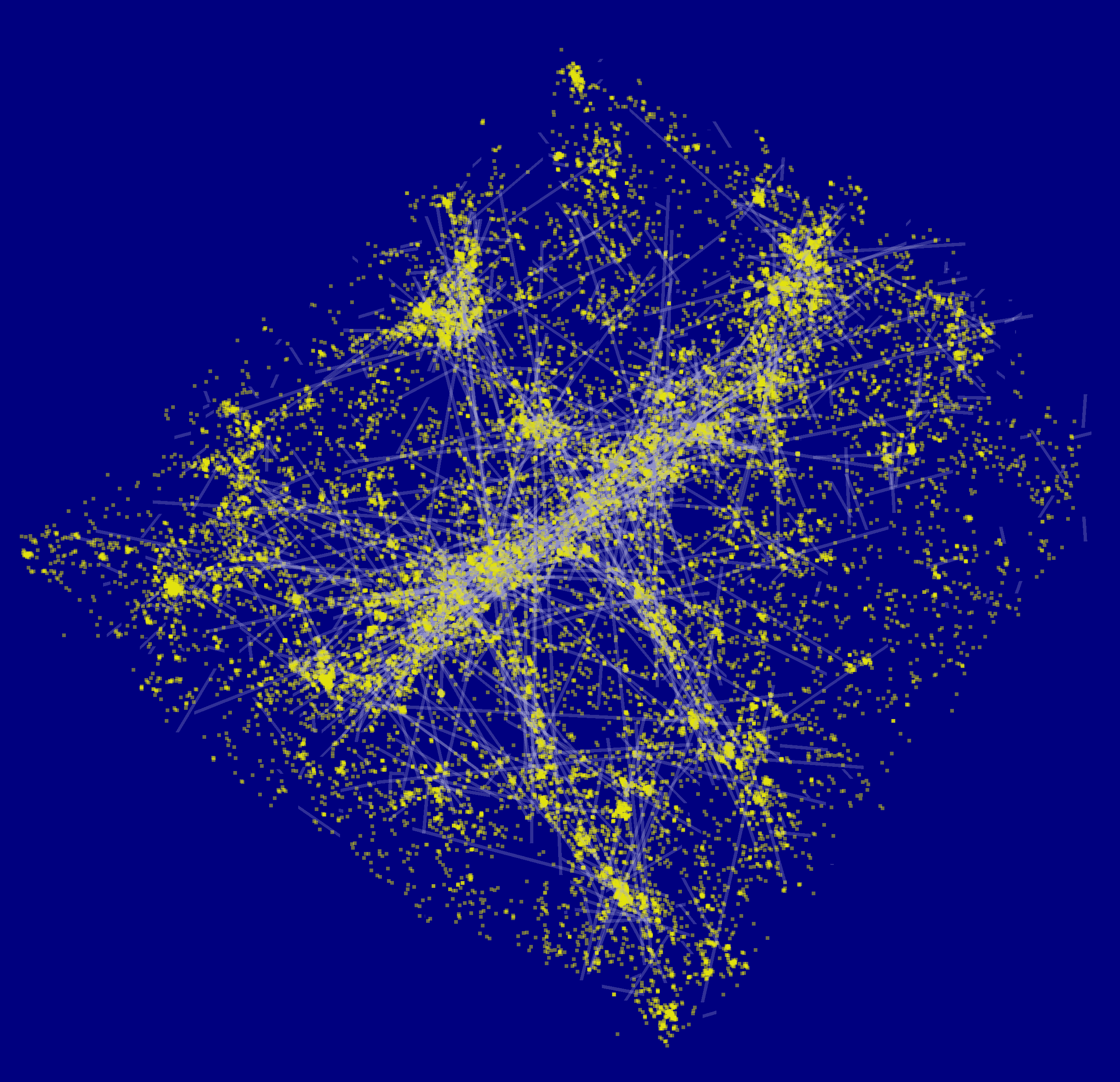}
 \includegraphics[width=.32\textwidth]{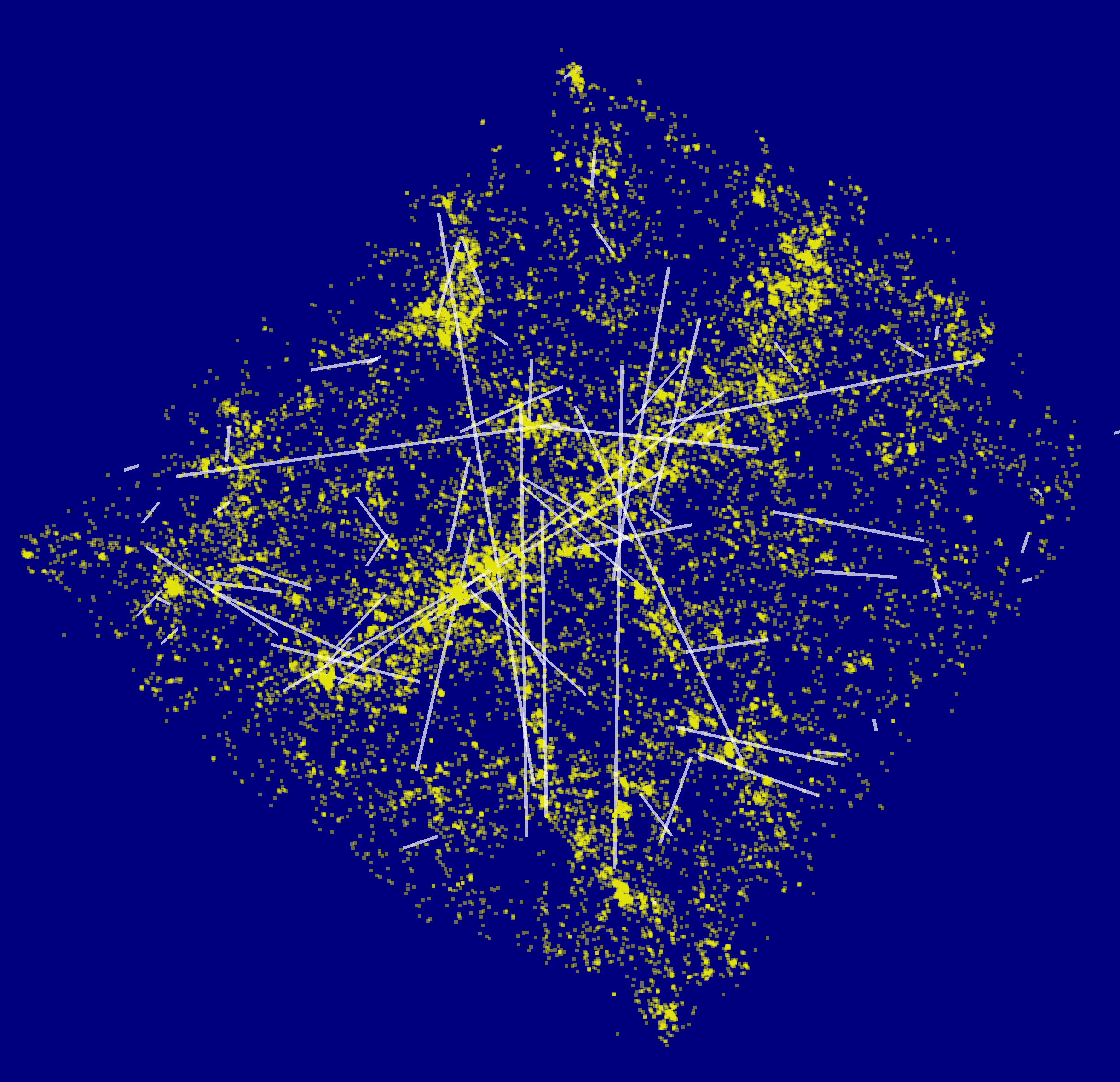}
 \includegraphics[width=.32\textwidth]{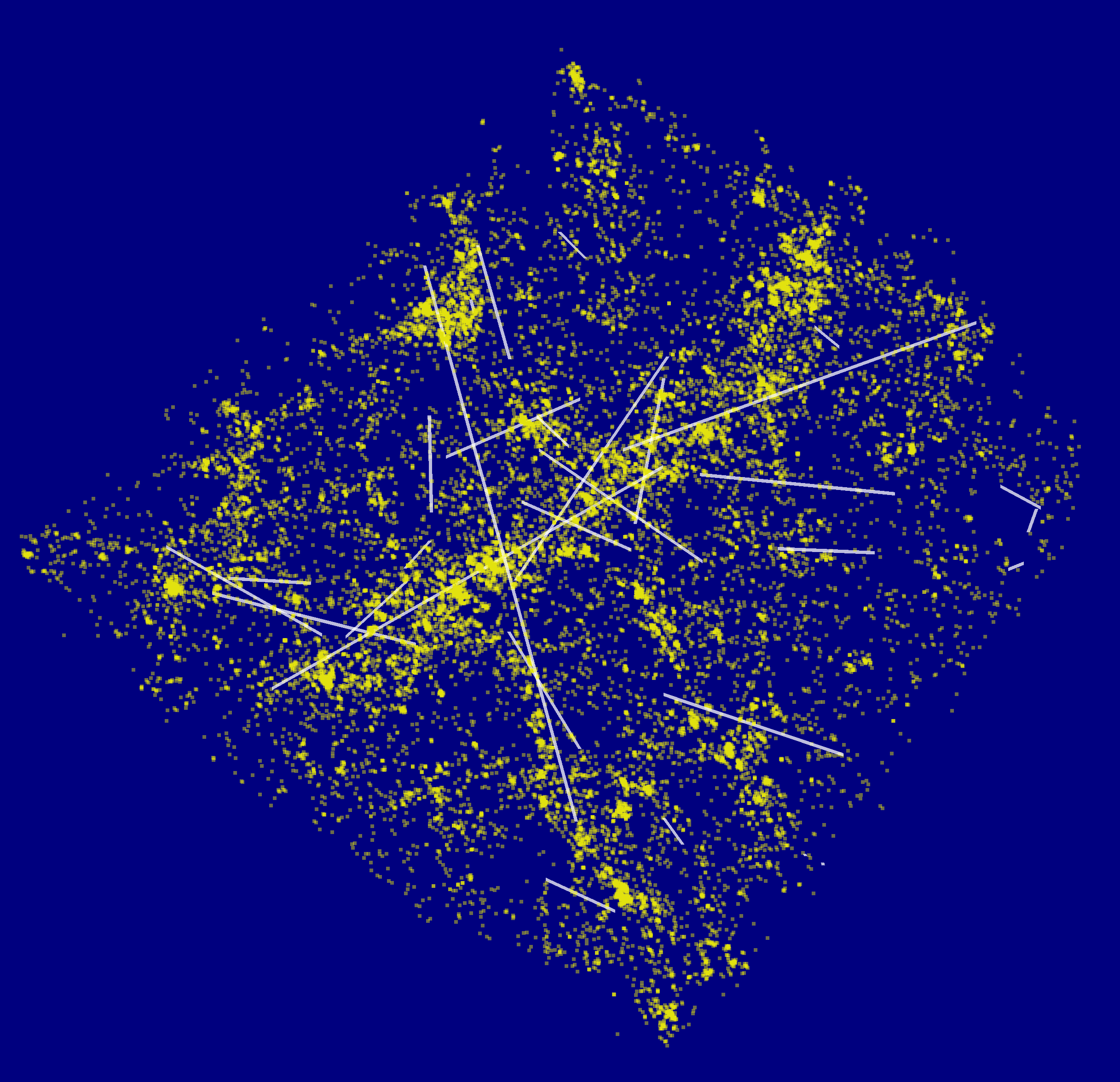}
 \caption{A thin slice of the global solution for locating Bisous filaments (white), and galaxies (yellow).}\label{fig::Stars}
\end{figure}

Galaxies are known to cluster along filaments and other low dimensional structures \cite{tempel2014galaxy,beygu2013interacting}. However mathematical extraction of these filaments is a challenging task. We might try to apply refinements of our algorithm to find those filaments. As a naive proof of concept, we have applied our method while setting $\rho_{2k+1} = 0$ for all $k$, so that the lines are disjoint. We use the data of \cite{tempel2016bisous} available at : \url{https://github.com/etempel/bisous}. The galaxies are represented by Dirac masses. Note that in our tests their masses $m$ are arbitrarily set to $\frac 1 n$ but the code enables other values. In Figure~\ref{fig::Stars}, computations are performed for $n=180\textrm{K}$ galaxies and a decreasing number of filaments $8\mathrm{K}$,$1\mathrm{K}$,$500$.
\clearpage

\bibliography{biblio}

\begin{thebibliography}{10}

\bibitem{aurenhammer1987power}
{\sc Aurenhammer, F.}
\newblock Power diagrams: properties, algorithms and applications.
\newblock {\em SIAM Journal on Computing 16}, 1 (1987), 78--96.

\bibitem{benamou2000computational}
{\sc Benamou, J.-D., and Brenier, Y.}
\newblock A computational fluid mechanics solution to the monge-kantorovich
  mass transfer problem.
\newblock {\em Numerische Mathematik 84}, 3 (2000), 375--393.

\bibitem{benamou2015iterative}
{\sc Benamou, J.-D., Carlier, G., Cuturi, M., Nenna, L., and Peyr{\'e}, G.}
\newblock Iterative bregman projections for regularized transportation
  problems.
\newblock {\em SIAM Journal on Scientific Computing 37}, 2 (2015),
  A1111--A1138.

\bibitem{bertsekas1999nonlinear}
{\sc Bertsekas, D.~P.}
\newblock {\em Nonlinear programming}.
\newblock Athena scientific Belmont, 1999.

\bibitem{beygu2013interacting}
{\sc Beygu, B., Kreckel, K., van~de Weygaert, R., van~der Hulst, J., and
  Van~Gorkom, J.}
\newblock An interacting galaxy system along a filament in a void.
\newblock {\em The Astronomical Journal 145}, 5 (2013), 120.

\bibitem{cgal}
\textsc{Cgal}, {C}omputational {G}eometry {A}lgorithms {L}ibrary.
\newblock http://www.cgal.org.

\bibitem{cuturi2013sinkhorn}
{\sc Cuturi, M.}
\newblock Sinkhorn distances: Lightspeed computation of optimal transport.
\newblock In {\em Advances in neural information processing systems\/} (2013),
  pp.~2292--2300.

\bibitem{de2012blue}
{\sc De~Goes, F., Breeden, K., Ostromoukhov, V., and Desbrun, M.}
\newblock Blue noise through optimal transport.
\newblock {\em ACM Transactions on Graphics (TOG) 31}, 6 (2012), 171.

\bibitem{Fred2018differentiation}
{\sc de~Gournay, F., Kahn, J., and Lebrat, L.}
\newblock Differentiation and regularity of semi-discrete optimal transport
  with respect to the parameters of the discrete measure.
\newblock {\em arXiv preprint arXiv:1803.00827\/} (2018).

\bibitem{Fred2017Gretsi}
{\sc de~Gournay, F., Kahn, J., Lebrat, L., and Pierre, W.}
\newblock Approches variationnelles pour le stippling: distances l2 ou
  transport optimal ?
\newblock In {\em GRETSI 2017 XXVI\/} (Sept. 2017).

\bibitem{Fred2018Projections}
{\sc de~Gournay, F., Kahn, J., Lebrat, L., and Pierre, W.}
\newblock Optimal transport approximation of measures.
\newblock {\em arXiv preprint arXiv:1804.08356\/} (2018).

\bibitem{dwyer1991higher}
{\sc Dwyer, R.~A.}
\newblock Higher-dimensional voronoi diagrams in linear expected time.
\newblock {\em Discrete \& Computational Geometry 6}, 3 (1991), 343--367.

\bibitem{fletcher2005barzilai}
{\sc Fletcher, R.}
\newblock On the barzilai-borwein method.
\newblock In {\em Optimization and control with applications}. Springer, 2005,
  pp.~235--256.

\bibitem{glowinski2014alternating}
{\sc Glowinski, R.}
\newblock On alternating direction methods of multipliers: a historical
  perspective.
\newblock In {\em Modeling, simulation and optimization for science and
  technology}. Springer, 2014, pp.~59--82.

\bibitem{hiller2003beyond}
{\sc Hiller, S., Hellwig, H., and Deussen, O.}
\newblock Beyond stippling—methods for distributing objects on the plane.
\newblock In {\em Computer Graphics Forum\/} (2003), vol.~22, Wiley Online
  Library, pp.~515--522.

\bibitem{hoff1999fast}
{\sc Hoff~III, K.~E., Keyser, J., Lin, M., Manocha, D., and Culver, T.}
\newblock Fast computation of generalized voronoi diagrams using graphics
  hardware.
\newblock In {\em Proceedings of the 26th annual conference on Computer
  graphics and interactive techniques\/} (1999), ACM Press/Addison-Wesley
  Publishing Co., pp.~277--286.

\bibitem{kitagawa2016convergence}
{\sc Kitagawa, J., M{\'e}rigot, Q., and Thibert, B.}
\newblock Convergence of a newton algorithm for semi-discrete optimal
  transport.
\newblock {\em arXiv preprint arXiv:1603.05579\/} (2016).

\bibitem{levy2015numerical}
{\sc L{\'e}vy, B.}
\newblock A numerical algorithm for l2 semi-discrete optimal transport in 3d.
\newblock {\em ESAIM: Mathematical Modelling and Numerical Analysis 49}, 6
  (2015), 1693--1715.

\bibitem{merigot2011multiscale}
{\sc M{\'e}rigot, Q.}
\newblock A multiscale approach to optimal transport.
\newblock In {\em Computer Graphics Forum\/} (2011), vol.~30, Wiley Online
  Library, pp.~1583--1592.

\bibitem{nesterov2013introductory}
{\sc Nesterov, Y.}
\newblock {\em Introductory lectures on convex optimization: A basic course},
  vol.~87.
\newblock Springer Science \& Business Media, 2013.

\bibitem{polyak1987introduction}
{\sc Polyak, B.~T.}
\newblock Introduction to optimization. translations series in mathematics and
  engineering.
\newblock {\em Optimization Software\/} (1987).

\bibitem{raydan1997barzilai}
{\sc Raydan, M.}
\newblock The barzilai and borwein gradient method for the large scale
  unconstrained minimization problem.
\newblock {\em SIAM Journal on Optimization 7}, 1 (1997), 26--33.

\bibitem{tempel2014galaxy}
{\sc Tempel, E., Kipper, R., Saar, E., Bussov, M., Hektor, A., and Pelt, J.}
\newblock Galaxy filaments as pearl necklaces.
\newblock {\em Astronomy \& Astrophysics 572\/} (2014), A8.

\bibitem{tempel2016bisous}
{\sc Tempel, E., Stoica, R.~S., Kipper, R., and Saar, E.}
\newblock Bisous model—detecting filamentary patterns in point processes.
\newblock {\em Astronomy and Computing 16\/} (2016), 17--25.

\bibitem{villani2008optimal}
{\sc Villani, C.}
\newblock {\em Optimal transport: old and new}, vol.~338.
\newblock Springer Science \& Business Media, 2008.

\end{thebibliography}
\bibliographystyle{acm}
\end{document}